\documentclass[11pt, leqno, a4paper,twoside]{amsart}
\usepackage[utf8]{inputenc}
\usepackage[numbers,sort&compress]{natbib}
\usepackage[draft]{changes}
\definechangesauthor[name=Paola, color=magenta]{PA}
\usepackage{amsmath}
\usepackage{siunitx}
\usepackage{amsfonts,bm}
\usepackage{dsfont}
\usepackage[colorlinks=true]{hyperref}
\usepackage{xcolor}
\usepackage{stmaryrd}
\usepackage{algorithm}
\usepackage{amssymb}
\usepackage{verbatim}
\usepackage{mathrsfs}
\usepackage{mathtools,amsmath}
\usepackage[margin=1.15in]{geometry}
\usepackage{tikz}
\usepackage{pgfplots}
\usepackage{capt-of}
\numberwithin{equation}{section}
\usepackage[noend]{algpseudocode}
\usepackage[capitalise, noabbrev]{cleveref}
\usepackage{caption}
\newtheorem{remark}{Remark}
\newtheorem*{assumption}{Mesh assumptions}
\newtheorem{theorem}{Theorem}

\newtheorem{lemma}{Lemma}

\newcommand*{\llbrace}{\{\mskip-5mu\{}
\newcommand*{\rrbrace}{\}\mskip-5mu\}}

\newcommand{\vecc}[1]{\mathbf{\underline{#1}}}
\newcommand{\norm}[1]{\left\lVert#1\right\rVert}
\newcommand{\seminorm}[1]{| #1|}
\definecolor{gray}{rgb}{0.5 0.5 0.5}
\definecolor{sgreen}{HTML}{40A441}
\definecolor{syellow}{HTML}{DFDE62}
\definecolor{spink}{HTML}{6D2A6F}
\definecolor{sred}{HTML}{B84041}
\definecolor{sblue}{HTML}{4040B2}
\definecolor{sgray}{HTML}{B6B6B6}

\title{A high-order discontinuous Galerkin method for nonlinear sound waves}
\author[P. F. Antonietti, I. Mazzieri, M. Muhr, V. Nikoli\'c, and B. Wohlmuth]{Paola. F. Antonietti$^1$, Ilario Mazzieri$^1$, Markus Muhr$^{\ast,2}$,\\ Vanja Nikoli\'c$^3$, and Barbara Wohlmuth$^2$}
\email{\href{paola.antonietti@polimi.it}{paola.antonietti@polimi.it}}
\email{\href{ilario.mazzieri@polimi.it}{ilario.mazzieri@polimi.it}}
\email{\href{mailto:muhr@ma.tum.de}{muhr@ma.tum.de}}
\email{\href{mailto:vanja.nikolic@ru.nl}{vanja.nikolic@ru.nl}}
\email{\href{mailto:wohlmuth@ma.tum.de}{wohlmuth@ma.tum.de}}

\thanks{$^*$Corresponding author: Markus Muhr, \href{mailto:muhr@ma.tum.de}{muhr@ma.tum.de}}
\begin{document}
\maketitle
\vspace*{-4mm}
\begin{center}
{\footnotesize
   $^1$MOX,  Dipartimento di Matematica, Politecnico di Milano, Milano, Italy  \\
   $^2$Department of Mathematics, Technical University of Munich, Germany \\
   $^3$Department of Mathematics, Radboud University, The Netherlands
}
\end{center}
\vspace{8mm}
\begin{abstract}
We propose a high-order discontinuous Galerkin scheme for nonlinear acoustic waves on polytopic meshes. To model sound propagation with and without losses, we use Westervelt's nonlinear wave equation with and without strong damping. Challenges in the numerical analysis lie in handling the nonlinearity in the model, which involves the derivatives in time of the acoustic velocity potential, and in preventing the equation from degenerating. We rely in our approach on the Banach fixed-point theorem combined with a stability and convergence analysis of a linear wave equation with a variable coefficient in front of the second time derivative. By doing so, we derive an \emph{a priori} error estimate for Westervelt's equation in a suitable energy norm for the polynomial degree $p \geq 2$. Numerical experiments carried out in two-dimensional settings illustrate the theoretical convergence results. In addition, we demonstrate efficiency of the method in a three-dimensional domain with varying medium parameters, where we use the discontinuous Galerkin approach in a hybrid way. 
\end{abstract}

\section{Introduction}
\indent Nonlinear sound waves arise in many different applications, such as medical ultrasound~\cite{duck2002nonlinear, rosnitskiy2015effect, maresca2017nonlinear}, fatigue crack detection~\cite{ryles2008comparative, sohn2014nonlinear}, or musical acoustics of brass instruments~\cite{campbell1999nonlinear, gilbert2007differences, myers2012effects}. Although considerable work has been devoted to their analytical studies~\cite{kaltenbacher2009global, kaltenbacher2011well, kawashima1992global, meyer2011optimal} and their computational treatment~\cite{hoffelner2001finite, kelly2018linear, resch2014two, tsuchiya1992simulation}, rigorous numerical analysis of nonlinear acoustic phenomena is still largely missing from the literature. The goal of our work is to develop a high-order discontinuous Galerkin  (DG) scheme for nonlinear sound waves that is rigorously justified through a stability and convergence analysis.\\
\indent The DG method was first introduced in the seventies for the numerical approximation of hyperbolic problems \cite{ReedHill73}, and, independently, in the context of elliptic~\cite{DouglasDupont76} and parabolic~\cite{Arnold82} equations. Since then DG methods have been successfully developed and applied to a wide range of problems arising in computational sciences and engineering; 
cf.\ the books \cite{hesthaven2007nodal,Riviere,di2011mathematical} for a comprehensive overview. In relation to our setting, we point out in particular the works on the Euler and Navier--Stokes equations~\cite{baumann1999discontinuous} and on a class of nonlinear elliptic and second-order hyperbolic problems~\cite{ortner2007discontinuous}. \\
\indent The finite-dimensional DG space consists of  piecewise discontinuous  polynomial functions defined over a computational tessellation of the domain. As a consequence, the DG paradigm can naturally support finite element spaces built upon  meshes consisting of arbitrarily shaped polygonal/polyhedral elements, thus generalizing the paradigm that stands at the basis of classical Finite Elements on triangles, quadrilaterals, or their combinations in two dimensions (2D), and tetrahedra, prisms, pyramids, and hexahedra or their combinations in three dimensions (3D), and gaining flexibility in the process of mesh generation. DG methods on polygonal/polyhedral grids (PolyDG methods for short) have received a lot of attention in the last years; we give here only an incomplete list~\cite{BaBoCoRe2012,BaBoCoSu2014,CaGeHo2014,CaDoGeHo2016,CaDoGeHo2017,AnCaCoDoGeGiHo2016} and refer the reader to the references therein for a comprehensive overview. In PolyDG methods, high order accuracy can be achieved in any space dimension by introducing suitable modal basis functions defined directly in the physical frame configuration. Finally, PolyDG methods can be seen as extensions of the classical DG approach and they are naturally oriented towards 3D scalable implementations.\\ 
\indent We organize the rest of the paper as follows. In Section~\ref{Sec:Continuous_problem}, we first discuss the continuous initial-boundary value problem for a classical model of nonlinear acoustics--Westervelt's wave equation. Section~\ref{Sec:Preliminaries} contains some theoretical preliminaries that are useful for the numerical analysis. In Section~\ref{Sec:Semi-discrete}, we propose and discuss a high-order discontinuous Galerkin scheme for the Westervelt equation. Section~\ref{Sec:Semi-discrete_Linear_Stability} is devoted to the stability analysis of a linearized semi-discrete problem and Section~\ref{Sec:Semi-discrete_Linear_Error} to its \emph{a priori} error analysis. In Section~\ref{Sec:Semi-discrete_Nonlinear_Error}, we use the Banach fixed-point theorem to prove an \emph{a priori} estimate for the approximate solution of the Westervelt equation. Section~\ref{Sec:Numerical_treatment} describes in detail our numerical solver. Finally, in Section~\ref{Sec:Numerical_results}, we carry out several numerical experiments, both in two and three dimensions, to illustrate the theory from previous sections. In a three-dimensional setting, we use the discontinuous Galerkin approach in a hybrid way to handle varying medium parameters. 
\section{The continuous problem}\label{Sec:Continuous_problem}
 We employ Westervelt's wave equation~\cite{westervelt1963parametric} to model nonlinear sound propagation,  given in terms of the acoustic velocity potential $\psi$ by
\begin{equation} \label{Westervelt_potential}
\begin{aligned}
(1-2k \dot{\psi})\ddot{\psi}-c^2 \Delta \psi - b \Delta \dot{\psi}=0.
\end{aligned}
\end{equation}
The constant $c$ denotes the speed of sound and $b$ is the so-called sound diffusivity. The constant $k$ is given by $k= \beta_a /c^2$, where $\beta_a$ is the coefficient of nonlinearity of the medium. For the derivation of nonlinear acoustic models and their physical background, we refer the interested reader to, e.g.,~\cite{Crighton, enflo2006theory, hamilton1998nonlinear}. After computing the acoustic velocity potential, the acoustic pressure $u$ can be obtained in a post-processing step via the relation $u=\varrho \dot{\psi}$, where $\varrho$ denotes the mass density of the medium. \\
\indent Westervelt's equation is a nonlinear acoustic wave equation, which we couple with initial conditions and homogeneous Dirichlet data, and investigate the following problem:
\begin{align} \label{Westervelt_ibvp_cont}
\begin{cases}
(1-2k \dot{\psi})\ddot{\psi}-c^2 \Delta \psi - b \Delta \dot{\psi}=0 \quad \text{in } \Omega \times (0,T],
\smallskip\\ \psi=0 \quad \text{ on }  \partial \Omega \times [0,T],  
\smallskip\\ (\psi, \dot{\psi})=(\psi_0, \psi_1) \quad \text{ on }   \Omega \times \{t=0\}
\end{cases}
\end{align}
on a bounded domain $\Omega \subset \mathbb{R}^d$ for $d \in \{2,3\}$ and for a given final time $T>0$. \\
\indent If $b>0$, then Westervelt's equation is strongly damped. With enough dissipation (i.e., $b$ large enough), it exhibits a parabolic-like behavior. The initial-boundary value problem \eqref{Westervelt_ibvp_cont} is then known to be globally well-posed for sufficiently small and smooth initial data on regular domains, provided that appropriate compatibility conditions at the initial time are satisfied. We refer to~\cite[Theorem 2.2]{kawashima1992global}, from which global well-posedness of \eqref{Westervelt_ibvp_cont} follows as a special case. We mention also the local-in-time well-posedness result from~\cite[Section 7]{kaltenbacher2014efficient} that relaxes the regularity assumptions on the initial data. \\
\indent If we consider propagation in inviscid media, then $b=0$ in \eqref{Westervelt_ibvp_cont}. It is expected and numerically observed~\cite{christov2007modeling, MKaltenbacher} that now smooth solutions of \eqref{Westervelt_ibvp_cont} exist only for a short time before the shock develops due to nonlinear steepening. A rigorous proof of the short-term well-posedness is available for propagation in unbounded domains as a particular case of a general quasi-linear hyperbolic system of second order in~\cite[Theorem 1]{hughes1977well}. For the inviscid Westervelt equation reformulated in terms of the acoustic pressure $u$, the local well-posedness on bounded domains follows from a special case of a general quasi-linear wave equation studied in~\cite[Theorem 4.1]{dorfler2016local}.\\
\indent In our numerical analysis, we intend to analyze both cases and assume that $b$ is non-negative. The analysis is valid as long as a sufficiently smooth solution of the original problem exists, i.e., up to the (possible) shock formation. Because we employ an energy method in the analysis, we have a delicate task of ensuring that all estimates we derive remain valid also in the absence of the strong damping, i.e., when $b=0$. \\
\indent We point out here another important feature of Westervelt's equation. The factor $1-2k \dot{\psi}$ in front of the second time derivative can degenerate if the acoustic pressure is too high. To avoid that this happens, we have to prove that $\dot{\psi}$ stays below $1/(2k)$. In the continuous analysis, this is commonly achieved by having sufficiently smooth data such that the solution space for the pressure embeds continuously into $L^\infty(\Omega)$ almost everywhere in time and by additionally assuming that the data are sufficiently small in an appropriate norm; see~\cite{kaltenbacher2014efficient, kaltenbacher2009global, kaltenbacher2011well, meyer2011optimal}. Our non-conforming discretization approach prevents this strategy. Since our approximate solution is only piecewise smooth, we have to rely on an inverse inequality to avoid degeneracy. On the other hand, we do not want a bound that degenerates as $h$ converges to $0$, and so we will need to involve additionally the (local) interpolant in the estimate and employ its approximation and stability properties. 
\section{Assumptions and preliminaries} \label{Sec:Preliminaries}
Let $\Omega \subset \mathbb{R}^d$ for $d \in \{2,3\}$ be a convex polygonal or polyhedral domain. We consider a family of meshes $\mathcal{T}_h$ made of disjoint open \emph{polygonal/polyhedral} elements $\kappa$ with diameter $h_{\kappa}$. \\
\indent Following~\cite{AnCaCoDoGeGiHo2016, CaDoGeHo2016, CaGeHo2014}, we introduce the concept of mesh \emph{interface}, defined as the intersection of the $(d-1)$-dimensional facets of two neighboring elements. When $d=3$, each interface consists of a general polygon which we assume can be decomposed into a set of co-planar triangles. We assume that a sub-triangulation of each interface is provided and we denote the set of all these triangles by $\mathcal{F}_h$. We then use the terminology \emph{face} to refer to one of the triangular elements in $\mathcal{F}_h$. When $d=2$, each interface simply consists of a line segment, so that the concept of faces and interfaces coincides in this case. We denote by $\mathcal{F}_h$ the set of all faces of $\mathcal{T}_h$, decomposed into the internal faces $\mathcal{F}_h^i$ and the boundary faces $\mathcal{F}_h^b$ so that $\mathcal{F}_h=\mathcal{F}^i_h \cup \mathcal{F}_h^b$.\\
\indent We assume a fixed uniform polynomial degree \(p \geq 1\) and introduce the following finite-dimensional space:
\[V_h=\{\psi  \in L^2(\Omega): \psi_{ \vert \kappa} \in \mathcal{P}_{p}(\kappa) \ \forall \kappa \in \mathcal{T}_h\},\]
where \(\mathcal{P}_{p}(\kappa) \) is the space of polynomials if total degree $p$ defined on $\kappa$,
as well as the broken Sobolev spaces \[H^{n}(\mathcal{T}_h)=\{\psi \in L^2(\Omega): \ \psi_{\vert \kappa} \in H^n(\kappa) \,\, \forall \kappa \in \mathcal{T}_h\}\] for $n \geq 1$. It is natural to employ the broken gradient operator $\nabla_h \cdot$ on the space $H^{1}(\mathcal{T}_h)$; see~\cite[Definition 1.21]{di2011mathematical}.\\
 \indent For sufficiently smooth $\psi$, we introduce jumps and averages on an interior face $F \in \mathcal{F}^i_h$, $F \subset \partial \kappa^+ \cap \partial \kappa ^-$ with $\kappa^+$ and $\kappa^-$ any two neighboring elements in $\mathcal{T}_h$, as follows:
\begin{align}\label{eq:def_jump_average}
\llbracket \psi \rrbracket =\psi^+ \mathbf{n}^+ + \psi^- \mathbf{n}^- , \qquad
\llbrace \psi \rrbrace=\frac{\psi^++\psi^-}{2},
\end{align}
where $\psi^{\pm}$ denotes the trace of $\psi$ on $F$ taken within the interior of $\kappa^\pm$, and $\mathbf{n}^{\pm}$ denotes the unit normal vector to $\partial \kappa^\pm$ pointing outwards from $\partial \kappa^\pm$. On the boundary face $F \in \mathcal{F}^b$, we set $\llbracket \psi \rrbracket =\psi \mathbf{n}$ and $\llbrace \psi \rrbrace = \psi$. For a (smooth enough) vector-valued function $\bm{\psi}$, definition \eqref{eq:def_jump_average} extends analogously.\\
\indent For later use, we also define here the stabilization function $\chi \in L^\infty(\mathcal{F}_h)$ as follows:
\begin{equation} \label{stabilization}
\begin{aligned}
\chi_{\vert F}= \begin{cases}
c^2 \, \beta \, \displaystyle \max_{\kappa \in \{\kappa^+, \kappa^-\}} \frac{p^2}{h_\kappa} \quad \quad &\text{for all  } F \in \mathcal{F}_h^i, \quad F \subset \partial \kappa^+ \cap \partial \kappa^-, \\[5mm]
\ c^2 \, \beta \, \dfrac{p^2}{h_\kappa}  &\text{for all  } F \in \mathcal{F}_h^b, \quad F \subset \partial \kappa,
\end{cases}
\end{aligned}
\end{equation}
The parameter $\beta>0$ will be chosen in a convenient manner in the following proofs.\\
\indent  For an open subset \(D\) of \(\mathbb{R}\), where \(d=1 ,2 ,3\), and a function $v \in H^n(D)$, where $n \geq 0$, we denote by \(\|v\|_{H^s(D)}\) and \(|v|_{H^s(D)}\) the standard norm and seminorm, respectively, with the convention that \(H^0(D)\equiv L^2(D)\). When $D \equiv \Omega$, we simply write $\|\nabla v\|_{H^n}$ and $|v|_{H^n}$. We use the short-hand notation
\[
\langle \psi, v \rangle_{\mathcal{F}}= \displaystyle \sum_{F \in \mathcal{F}} (\psi, v)_{L^2(F)}, \quad \|\psi\|_{\mathcal{F}}=\langle \psi, \psi \rangle^{1/2}_{\mathcal{F}} 
\]
for a generic collection of faces $\mathcal{F} \subset \mathcal{F}_h$, and regular enough functions $\psi$ and $v$. Here \((\cdot, \cdot)_{L^2(F)}\) denotes the inner product in $L^2(F)$.\\
\indent We occasionally use the notation $x \lesssim y$ and $x \gtrsim y$ instead of $x \leq Cy$ and $x \geq Cy$, respectively, when the hidden constant $C>0$ does not depend on the coefficients in the equation $c$, $b$, and $k$, the mesh size, and the number of faces of a mesh element, but can depend on the polynomial degree $p$ and the final time $T$.

\subsection{Grid assumptions and preliminary estimates}
Throughout the paper,  we make the following assumptions on the family of polytopic decompositions $\mathcal{T}_h$, which allow to extend the trace-inverse  and inverse inequalities on simplices to polytopic elements.
\begin{assumption}\label{ass1}
For any $\kappa \in \mathcal{T}_h$, we assume that \[ h_\kappa^d \geq |\kappa| \gtrsim h_\kappa^d\] for $d=2, 3$,
where $|\kappa|$ denotes the Hausdorff measure of $\kappa \in \mathcal{T}_h$. We also assume that there exists a positive number $m$, such that every polytopic element $\kappa \in \mathcal{T}_h$ admits a sub-triangulation into at most $\mathtt{m}_{\kappa} \leq  m$ shape-regular simplices $\mathtt{s}_i$ for $i=1,2,\ldots, \mathtt{m}_{\kappa}$, such that \[\bar{\kappa}=\cup_{i=1}^{\mathtt{m}_{\kappa}} \bar{\mathtt{s}}_i \ \text{  and  } \ |\mathtt{s}_i| \gtrsim |\kappa|, \] where the hidden constant is independent of $\kappa$ and $\mathcal{T}_h$. Finally, we assume that 
\begin{align*} 
\displaystyle \frac{\max_{\kappa} h_\kappa}{\min_{\kappa} h_{\kappa}} \lesssim 1.
\end{align*}
\end{assumption}
Under these mesh assumptions, the following trace-inverse and inverse inequalities hold on polytopic domains.
\begin{lemma}\label{lem:inversepoly}
For any $v \in \mathcal{P}_{p}(\kappa)$, $\kappa \in \mathcal{T}_h$, the following trace-inverse and inverse inequalities hold:
\begin{align}
\norm{v}_{L^2(\partial \kappa)} &\lesssim h_{\kappa}^{-1/2} \, \norm{v}_{L^2(\kappa)}, \label{eq:tracepoly}\\
\norm{v}_{L^{\infty}(\kappa)} &\lesssim h_{\kappa}^{-d/2} \, \norm{v}_{L^2(\kappa)}. \label{eq:inversepoly}
\end{align}
\end{lemma}
\begin{proof}
The statement follows from, e.g.,~\cite[Lemma 6]{CaDoGeHo2017} combined with our mesh assumptions.
\end{proof}

\begin{remark}[On the mesh assumptions]\label{rem:remark_trace_inverse_poly} We choose to simplify our mesh assumptions for the clarity of exposition. However, for the trace-inverse inequality \eqref{eq:tracepoly} to hold, these assumptions are slightly more restrictive than needed, and can be weakened by employing the arguments of~\cite{CaDoGe2017}. Indeed,  the inequality  holds provided that, for any $\kappa \in \mathcal{T}_h$, there exists a set of non-overlapping $d$-dimensional simplices $\kappa_{\flat}^F \subset \kappa$ such that, for any face $F \subset \partial \kappa$, \(\overline{F} = \partial \overline{\kappa} \cap \partial \overline{\kappa}_{\flat}^F\), and $\bigcup_{F \subset \partial \kappa} \overline{\kappa_{\flat}^F} \subset \overline{\kappa}$, and the diameter $h_{\kappa}$ of $\kappa$ can be bounded by
\[ h_{\kappa} \lesssim \frac{d |\kappa_{\flat}^F |}{|F|}\] for all $F \subset \partial \kappa$,
where $|F|$ and $|\kappa_{\flat}^F|$ denote the Hausdorff measure of $F$ and $\kappa_{\flat}^F$, respectively. This latter assumption does not put a restriction on either the number of faces that an element possesses, or indeed the measure of a face of an element $\kappa \in \mathcal{T}_h$, relative to the measure of the element itself; cf.\ also
\cite{AnCaCoDoGeGiHo2016, CaDoGeHo2016,CaDoGeHo2017, CaGeHo2014}. As pointed out in \cite{CaDoGe2017}, meshes obtained by agglomeration of a finite number of polygons that are uniformly star-shaped with respect to the largest inscribed ball will automatically satisfy the above weak requirement.\\
\indent The inverse inequality \eqref{eq:inversepoly} also holds under weaker assumptions: if, for any point $x \in \kappa$, there exists a shape-regular simplex  containing $x$  and contained in $\kappa$, with diameter comparable to that of $\kappa$. In other words, for any point $x \in \kappa$, there exists $\mathtt{s}_{\kappa}(x)$, such that $x \in \mathtt{s}_{\kappa}(x) \subseteq \kappa$ and $h_{\mathtt{s}(x)} \gtrsim h_\kappa$. The proof follows by relying on the $L^{\infty}$ inverse estimates on shape-regular simplexes; cf.~\cite[equation (3.8)]{Georgoulis_08},  \cite[Theorem 4.76]{Schwab98}.
\end{remark}
\subsection{Interpolation bounds on polytopic meshes} For future reference, we also state here the specific interpolation bounds on polytopic meshes we will rely on in the proofs.
\begin{lemma}\label{lem:interpolation_bounds}
		Let  $v \in H^{n}(\kappa)$, where $\kappa \in \mathcal{T}_h$. Then, there exists $\Pi_{\kappa,p}: H^{n}(\kappa) \rightarrow \mathcal{P}_p(\kappa)$ such that
		\begin{equation}\label{eq:interpolation_bounds}
		\begin{aligned}
	\norm{v-\Pi_{\kappa, p} v}_{L^2(\kappa)}& \lesssim&& h_{\kappa}^{\mu}
		\seminorm{v}_{H^{n}(\kappa)}, \ &&&n\geq 0, \\
			\seminorm{v-\Pi_{\kappa, p} v}_{H^{1}(\kappa)}& \lesssim&& h_{\kappa}^{\mu-1}
		\seminorm{v}_{H^{s}(\kappa)}, \ &&&n\geq 1, \\
\norm{v-\Pi_{\kappa, p} v}_{L^\infty(\kappa)}& \lesssim&&  h_{\kappa}^{\mu-d/2}  \seminorm{v}_{H^n(\kappa)}, \ &&&n>d/2,
		\end{aligned}
		\end{equation}
		where $\mu = \min \{n, p+1\}$. 
\end{lemma}
\begin{proof}
The statement follows by employing our mesh assumptions and classical interpolation bounds on quadrilateral/hexahedral  and simplicial elements; cf.~\cite{Babuska_Suri_1987b, CaGeHo2014, brenner2007mathematical}.
\end{proof}
We can now also state a result on the interpolation error for time-dependent piecewise smooth functions.  Let $\psi \in C([0,T]; H^n(\mathcal{T}_h))$, where $n \geq 2$. For any time $t \in [0,T]$, we define the global interpolant $\psi_I$ element-wise as 
\begin{align} \label{global_interpolant_definition}
{\psi_I}_{\vert \kappa}=\psi_{I,\kappa}, \quad \kappa \in \mathcal{T}_h,
\end{align}
 where $\psi_{I,\kappa}=\Pi_{\kappa, p}\psi$ is the local interpolant of~\Cref{lem:interpolation_bounds}.
\begin{lemma}\label{Lemma:Interpolant} 
Let $\psi \in C([0,T]; H^n(\mathcal{T}_h))$, where $n \geq 2$. Then, there exists an interpolant $\psi_I \in C([0,T]; V_h)$, defined as in \eqref{global_interpolant_definition}, such that the error $e_I=\psi-\psi_I$ satisfies
\begin{equation*}
\begin{aligned}
c^2\|\nabla_h e_I(t)\|^2_{L^2} +\|\sqrt{\chi} \llbracket e_I(t) \rrbracket \,\|^2_{\mathcal{F}_h}
  \lesssim \, c^2 \sum_{\kappa \in \mathcal{T}_h} h_\kappa^{2 \mu-2}\, \seminorm{ \psi(t)}^2_{H^n(\kappa)},
  \end{aligned}
  \end{equation*}
for all $t \in [0, T]$, where $\mu = \min \{n, p+1\}$ and \(\chi\) is defined in \eqref{stabilization}.  Moreover, the following estimate holds:
\begin{equation} \label{interpolant_error_2}
\begin{aligned}
 \|\chi^{-1/2} \llbrace \nabla_h e_I(t) \rrbrace \,\|^2_{\mathcal{F}_h} \lesssim \,\frac{1}{c^2} \sum_{\kappa \in \mathcal{T}_h} h_\kappa^{2 \mu-2}\, \seminorm{\psi(t)}^2_{H^n(\kappa)}, \quad \text{for }\ 0 \leq t \leq T.  \\
\end{aligned}
\end{equation}
\end{lemma}
\begin{proof}
The statement follows by relying on the mesh assumptions, estimates \eqref{eq:interpolation_bounds}, and the following multiplicative trace inequality on shape-regular simplices $\mathtt{s}$:
\begin{equation*}
\begin{aligned}
&
\norm{\eta}_{L^2(\partial \mathtt{s})}^2
\lesssim
\norm{\eta}_{L^2(\mathtt{s})}
\norm{\nabla \eta}_{L^2(\mathtt{s})}
+ h_\mathtt{s}^{-1}\norm{\eta}_{L^2(\mathtt{s})}^2
&& \ \text{for all } \eta \in H^{1}(\mathtt{s});
\end{aligned}
\end{equation*}
cf.~\cite{Babuska_Suri_1987b} and \cite[Lemma 33]{CaDoGeHo2017}.
\end{proof}
\section{The DG approximation in space of the Westervelt equation} \label{Sec:Semi-discrete}
\indent In this section, we introduce and discuss the semi-discrete approximation of the initial-boundary value problem \eqref{Westervelt_ibvp_cont} for the Westervelt equation. To motivate our approximate weak form, we rewrite the Westervelt equation in \eqref{Westervelt_ibvp_cont}as
\begin{equation*} 
\begin{aligned}
(1-2k \dot{\psi})\ddot{\psi}-c^2 \Delta (\psi +\tfrac{b}{c^2}\dot{\psi})=0.
\end{aligned}
\end{equation*}
Together with the fact that sound diffusivity $b$ is relatively small in realistic applications, this suggests to introduce an auxiliary state
\begin{equation} \label{psi_underline}
\begin{aligned}
\tilde{\psi}= \psi+\frac{b}{c^2}\dot{\psi},
\end{aligned}
\end{equation}
which allows us to formally write the Westervelt equation as
\begin{equation} \label{Westervelt_underline}
\begin{aligned}
(1-2k \dot{\psi})\ddot{\psi}-c^2 \Delta \tilde{\psi}=0.
\end{aligned}
\end{equation}
If wave propagates through inviscid media, then the auxiliary state $\tilde{\psi}$ is equal to the acoustic velocity potential $\psi$ and \eqref{Westervelt_underline} reduces to the inviscid Westervelt equation. The auxiliary state helps us to unify in our analysis propagation in inviscid and non-inviscid media and assume that $b \geq 0$.\\
\indent We are interested in the solutions of this problem in the sense of the equation 
\begin{equation} \label{Westervelt_ibvp_weak}
\begin{aligned}
((1-2k \dot{\psi}) \ddot{\psi}, v)_{L^2}+a(\tilde{\psi}, v) = 0 
\end{aligned}
\end{equation}
being satisfied for all $v \in H^1(\Omega)$ and all times $t \in (0,T]$, with $(\psi, \dot{\psi})\vert_{t=0}=(\psi_0, \psi_1)$. The bilinear form $a:H^1(\Omega) \times H^1(\Omega)\rightarrow \mathbb{R}$ is given by
\begin{equation*} \label{bilinear_form}
\begin{aligned}
a(\psi, v)=c^2 (\nabla \psi , \nabla v)_{L^2}. 
\end{aligned}
\end{equation*}
We introduce the corresponding DG bilinear form $a_h: H^1(\mathcal{T}_h) \times V_h \rightarrow \mathbb{R}$ by
\begin{equation*} \label{discrete_bilinear_form}
\begin{aligned}
a_h(\psi, v_h)=& \,\begin{multlined}[t]c^2 (\nabla_h \psi , \nabla_h v_h)_{L^2} 
- \langle \llbrace c^2 \nabla_h \psi \rrbrace, \llbracket v_h \rrbracket \rangle_{\mathcal{F}_h} \\
- \langle \llbracket \psi \rrbracket,  \llbrace c^2 \nabla_h v_h \rrbrace \rangle_{\mathcal{F}_h} \, +\langle \chi  \llbracket \psi \rrbracket ,  \llbracket v_h \rrbracket  \rangle_{\mathcal{F}_h}, \end{multlined}
\end{aligned}
\end{equation*}
where the stabilization function $\chi$ is defined as in \eqref{stabilization}. We then look for the approximate solution $\psi_h \in C^2([0,T]; V_h)$ of equation \eqref{Westervelt_ibvp_weak}, such that 
\begin{equation} \label{Westervelt_ibvp_discr}
\begin{aligned}
 ((1-2k \dot{\psi}_h) \ddot{\psi}_h, v_h)_{L^2}+a_h(\tilde{\psi}_h, v_h) = 0 
\end{aligned}
\end{equation}
holds for all $ v_h \in V_h,\ 0 < t \leq T$, supplemented with the approximate initial data
$$(\psi_h(0), \dot{\psi}_{h}(0))=(\psi_{0, h}, \psi_{1, h}) \in V_h \times V_h.$$
In equation \eqref{Westervelt_ibvp_discr},  we have used, analogously to \eqref{psi_underline}, the notation
\begin{equation} \label{psi_h_underline}
\begin{aligned}
\tilde{\psi}_h= \psi_h+\frac{b}{c^2}\dot{\psi}_h,
\end{aligned}
\end{equation}
and therefore our weak form \eqref{Westervelt_ibvp_discr} is equivalent to
\begin{equation} \label{weak_form_written_out}
\begin{aligned}
\begin{multlined}[t]((1-2k \dot{\psi}_h) \ddot{\psi}_h, v_h)_{L^2}+c^2 (\nabla_h \psi_h , \nabla_h v_h)_{L^2} +b (\nabla_h \dot{\psi}_h , \nabla_h v_h)_{L^2} \\
- \langle \llbrace c^2 \nabla_h \psi_h \rrbrace, \llbracket v_h \rrbracket \rangle_{\mathcal{F}_h}- \langle \llbrace b \nabla_h \dot{\psi}_h \rrbrace, \llbracket v_h \rrbracket \rangle_{\mathcal{F}_h} \\
- \langle \llbracket \psi_h \rrbracket,  \llbrace c^2 \nabla_h v_h \rrbrace \rangle_{\mathcal{F}_h} - \langle \llbracket \tfrac{b}{c^2} \dot{\psi}_h \rrbracket,  \llbrace c^2 \nabla_h v_h \rrbrace \rangle_{\mathcal{F}_h} \\ +\langle \chi  \llbracket \psi_h \rrbracket ,  \llbracket v_h \rrbracket  \rangle_{\mathcal{F}_h}+ \langle \chi  \llbracket \tfrac{b}{c^2} \dot{\psi}_h \rrbracket ,  \llbracket  v_h \rrbracket  \rangle_{\mathcal{F}_h}=0. \end{multlined}
\end{aligned}
\end{equation}
Recall that in \eqref{stabilization} the stabilization function $\chi$ has a $c^2$ scaling, and thus the two stabilization terms in \eqref{weak_form_written_out} effectively scale by $c^2$ and $b$. \\
\indent We note that in the case of sound propagation without losses, where $b=0$ in the Westervelt equation, \eqref{Westervelt_ibvp_discr} corresponds to the standard DG formulations for second-order undamped wave equations; see, for example,~\cite{grote2006discontinuous}.
\section{Analysis of the linearized semi-discrete problem} \label{Sec:Semi-discrete_Linear_Stability}
As a first step in the analysis, we consider a non-degenerate linearization of \eqref{Westervelt_ibvp_cont} that is given by the following initial-boundary value problem for a linear strongly damped wave equation:
\begin{align} \label{Westervelt_ibvp_cont_linearized}
\begin{cases}
\alpha(x,t)\ddot{\psi}-c^2 \Delta \tilde{\psi}=0 \quad \text{in } \Omega \times (0,T],
\smallskip\\ \psi=0 \quad \text{ on }  \partial \Omega \times [0,T],  
\smallskip\\ (\psi, \dot{\psi})=(\psi_0, \psi_1) \quad \text{ on }   \Omega \times \{t=0\},
\end{cases}
\end{align}
where it is assumed that there exist $\alpha_0$, $\alpha_1>0$ such that 
\begin{align*} \label{assumptions_alpha}
\alpha_0 \leq \alpha(x,t) \leq \alpha_1 \ \text{ in } \overline{\Omega} \times [0,T],
\end{align*}
and the relation \eqref{psi_underline} holds. Sufficient conditions for the well-posedness of \eqref{Westervelt_ibvp_cont_linearized} in the case that $b>0$ can be found in, e.g.,~\cite[Proposition 3.2]{kaltenbacher2014efficient}. The weak form of this problem is given by
\begin{equation*}  \label{Westervelt_lin_weak}
\begin{aligned}
(\alpha \ddot{\psi}, v)_{L^2}+a(\tilde{\psi}, v) = 0 
\end{aligned}
\end{equation*}
for all $v \in H^1(\Omega)$, $0<t\leq T$ with $(\psi, \dot{\psi})\vert_{t=0}=(\psi_0, \psi_1)$. We analyze its semi-discrete approximation, given by equation
\begin{equation}  \label{Westervelt_lin_ibvp_discr}
\begin{aligned}
 (\alpha_h \ddot{\psi}_h, v_h)_{L^2}+a_h(\tilde{\psi}_h, v_h) = 0 
\end{aligned}
\end{equation}
which should hold for all $ v_h \in V_h,\ 0 < t \leq T$, supplemented with the approximate initial data 
$$(\psi_h(0), \dot{\psi}_{h}(0))=(\psi_{0, h}, \psi_{1, h}) \in V_h \times V_h.$$ 
In equation \eqref{Westervelt_lin_ibvp_discr}, the coefficient $\alpha_h$ denotes a discrete version of the coefficient $\alpha$ such that 
\begin{align} \label{assumptions_alpha_h}
\alpha_0 \leq \alpha_h(x,t) \leq \alpha_1 \ \text{ in } \overline{\Omega} \times [0,T].
\end{align}
\indent The main idea behind studying this linearized problem is to later choose $$\alpha_h=1-2k \dot{w}_h$$ with $w_h$ in a neighborhood of $\psi$,
 and define a map $\mathcal{J}: w_h \mapsto \psi_h$, where $\psi_h$ solves the linear semi-discrete problem \eqref{Westervelt_lin_ibvp_discr}. The fixed point of this map will be the solution of the nonlinear problem \eqref{Westervelt_ibvp_discr}. Our approach here follows, in spirit, the strategy taken in~\cite{ortner2007discontinuous}, where nonlinear hyperbolic systems in divergence form are considered.
\subsection{Existence and stability}
\indent Our first task is to prove that the semi-discrete problem \eqref{Westervelt_lin_ibvp_discr} has a unique solution that remains bounded in a suitable energy norm. We begin by recalling a useful inequality for functions in $V_h$.
\begin{lemma} \label{Lemma:Ineq}
For any $v_h \in V_h$, the following inequality holds:
\begin{equation*} \label{aux_ineq_1}
\begin{aligned}
\|\chi^{-1/2} \llbrace \nabla_h v_h \rrbrace\|_{\mathcal{F}_h} \lesssim  \frac{1}{c\sqrt{\beta}}\|\nabla_h v_h\|_{L^2},
\end{aligned}
\end{equation*}
where $\beta>0$ is the stability parameter that appears in the definition \eqref{stabilization} of the stabilization function $\chi$. 
\end{lemma}
\begin{proof}
The statement follows by a straightforward modification of the arguments in  \cite[Lemma 3.2]{AnMa2018}; cf.\ also~\cite{AntoniettiAyusoMazzieriQuarteroni_2016}.
\end{proof}
By relying on Lemma~\ref{Lemma:Ineq}, we can show that
\begin{equation} \label{helpful_ineq}
\begin{aligned}
c^2 \|\langle \llbrace  \nabla_h v_h  \rrbrace , \llbracket v_h \rrbracket \rangle_{\mathcal{F}_h}\| \lesssim& \, c \frac{1}{\sqrt{\beta}}\|\nabla_h v_h \|_{L^2}\|\sqrt{\chi}\, \llbracket v_h \rrbracket \, \|_{\mathcal{F}_h} \\
\lesssim&\,   c^2\frac{1}{4 \varepsilon} \frac{1}{\beta}\|\nabla_h v_h \|^2_{L^2}+\varepsilon \|\sqrt{\chi}\, \llbracket v_h \rrbracket \, \|^2_{\mathcal{F}_h}
\end{aligned}
\end{equation}
for functions $v_h \in V_h$, which we will rely on in the upcoming proofs by choosing $\varepsilon>0$ in a convenient manner.\\
\indent In order to state our results, we introduce the discrete energy function
\begin{equation*} \label{energy}
\begin{aligned}
E[\psi_h](t):=&\,\begin{multlined}[t] \|\sqrt{\alpha_h(t)}\, \dot{\psi}_h(t)\|^2_{L^2}+\tfrac{b}{c^2}\int_0^t \|\sqrt{\alpha_h}\, \ddot{\psi}_h\|_{L^2}^2 \, \textup{d}s \\[1mm]
	+ c^2\|\nabla_h \tilde{\psi}_h(t)\|^2_{L^2}+\|\sqrt{\chi}\, \llbracket \tilde{\psi}_h(t) \rrbracket \,\|^2_{\mathcal{F}_h}, \end{multlined}
		\end{aligned}
\end{equation*}
for $t \in [0,T]$. We note that in the case $b=0$, we have $\tilde{\psi}_h=\psi_h$, and we recover the energy of undamped linear wave equations. 
\begin{theorem} \label{Thm:LinStability} Let $c>0$, $b \geq 0$, and let $T>0$ be a fixed time horizon. Let the coefficient $\alpha_h \in C^1([0,T]; V_h)$ be such that the non-degeneracy condition \eqref{assumptions_alpha_h} holds, where $\alpha_0$ and $\alpha_1$ are independent of the discretization parameters. Moreover, assume that there exists $\gamma \in (0,1)$ such that
\begin{equation} \label{alpha_t_small}
\left \|\dot{\alpha}_{h}/\alpha_h \right\|_{L^1 L^\infty}  \leq \gamma.
\end{equation}
	\indent Then the semi-discrete problem \eqref{Westervelt_lin_ibvp_discr} has a unique solution $\psi_h$ such that it holds
	\begin{equation} \label{WellpLin:final_est}
	\begin{aligned}
\max_{t \in [0,T]} E[\psi_h](t) \leq C_{\textup{Th}\ref{Thm:LinStability}} \,E[\psi_h](0),
	\end{aligned}
	\end{equation}	
provided that the parameter $\beta$ in \eqref{stabilization} is sufficiently large. The constant $C_{\textup{Th}\ref{Thm:LinStability}}>0$ does not depend on the mesh size, the number of faces of a mesh element, or the coefficients in the equation, but depends on the polynomial degree.
\end{theorem}
\begin{proof}
Because the problem is non-degenerate and $\alpha_h \in C^1([0,T]; V_h)$, local-in-time existence of a solution $\psi_h \in C^2([0, T_h]; V_h)$ for some $T_h \leq T$ follows by relying on the standard theory of linear ordinary differential equations; cf. \cite[Theorem 1.44]{roubivcek2013nonlinear}  and \cite[Theorem 4.2]{nikolic2019priori}. The upcoming energy estimate will allow us to extend the existence interval to $[0,T]$. \\
\indent We next focus on proving stability. In energy analysis of second-order wave equations, the first time derivative of the solution is a natural choice of test function. However, due to the presence of a varying coefficient $\alpha_h$ in our case, we would need to additionally test with a suitably scaled second time derivative.  We combine these ideas and choose $ v_h=\dot{\tilde{\psi}}_h= \dot{\psi}_h+\tfrac{b}{c^2}\ddot{\psi}_h$ as a test function.\\
\indent Taking $v_h= \dot{\tilde{\psi}}_h$ in \eqref{Westervelt_lin_ibvp_discr}, integrating over $(0,t)$, where $t \leq T_h$, and performing integration by parts with respect to time leads to the identity
\begin{equation} \label{LinWellp:Identity}
\begin{aligned}
&\int_0^t (\alpha_h \ddot{\psi}_h, \dot{\tilde{\psi}}_h)_{L^2}\, \textup{d}s+\tfrac{1}{2}c^2 \|\nabla_h \tilde{\psi}_h(s)\|^2_{L^2}\, \Bigl \vert_0^t +\tfrac12\|\sqrt{\chi}\, \llbracket \tilde{\psi}_h(s) \rrbracket\|^2_{\mathcal{F}_h}\Bigl \vert_0^t\\[1mm]
=& \, \langle \llbrace c^2 \nabla_h \tilde{\psi}_h(s)  \rrbrace , \llbracket \tilde{\psi}_h(s) \rrbracket \rangle_{\mathcal{F}_h}\, \Bigl \vert_0^t.
\end{aligned}
\end{equation} 
Above, we have made use of the fact that
\begin{equation*} \label{c_terms_1}
\begin{aligned}
\langle \llbrace c^2 \nabla_h \tilde{\psi}_h  \rrbrace , \llbracket \dot{\tilde{\psi}}_h \rrbracket \rangle_{\mathcal{F}_h}+\langle \llbracket \tilde{\psi}_h \rrbracket, \llbrace c^2 \nabla_h \dot{\tilde{\psi}}_h  \rrbrace  \rangle_{\mathcal{F}_h}=\frac{\textup{d}}{\textup{d}t}\, \langle \llbrace c^2 \nabla_h \tilde{\psi}_h  \rrbrace , \llbracket \tilde{\psi}_h \rrbracket \rangle_{\mathcal{F}_h}.
\end{aligned}
\end{equation*}
We can employ Lemma~\ref{Lemma:Ineq} and inequality \eqref{helpful_ineq} to obtain
\begin{equation*}  \label{c_terms_2}
\begin{aligned}
&\tfrac{1}{2}c^2 \|\nabla_h \tilde{\psi}_h\|^2_{L^2}+\tfrac12 \|\sqrt{\chi}\, \llbracket \tilde{\psi}_h \rrbracket \,\|^2_{\mathcal{F}_h}-c^2\langle \llbrace  \nabla_h \tilde{\psi}_h  \rrbrace , \llbracket \tilde{\psi}_h \rrbracket \rangle_{\mathcal{F}_h}\\
\geq & \, C_1( c^2 \|\nabla_h \tilde{\psi}_h\|^2_{L^2}+\|\sqrt{\chi}\, \llbracket \tilde{\psi}_h\rrbracket \,\|^2_{\mathcal{F}_h})
\end{aligned}
\end{equation*}
for all $t \in [0,T_h]$, provided that parameter $\beta$ in \eqref{stabilization} is sufficiently large. Similarly, we obtain
\begin{equation*}
\begin{aligned}
& \tfrac{1}{2}c^2 \|\nabla_h \tilde{\psi}_h(0)\|^2_{L^2}+\|\sqrt{\chi}\, \llbracket \tilde{\psi}_h(0) \rrbracket\, \|^2_{\mathcal{F}_h}- \langle \llbrace c^2 \nabla_h \tilde{\psi}_h(0)  \rrbrace , \llbracket \tilde{\psi}_h(0) \rrbracket \rangle_{\mathcal{F}_h}\\ \leq&\, C_2(c^2 \|\nabla_h \tilde{\psi}_h(0)\|^2_{L^2}+\|\sqrt{\chi}\, \llbracket \tilde{\psi}_h(0) \rrbracket \,\|^2_{\mathcal{F}_h}).
\end{aligned}
\end{equation*}
The constants $C_1$, $C_2>0$ above are independent of $c$, $b$, or the mesh size, but depend on the polynomial degree \(p\). It remains to estimate the $\alpha_h$-term in \eqref{LinWellp:Identity}. We recall how the auxiliary state $\tilde{\psi}_h$ is defined in \eqref{psi_h_underline} and employ integration by parts with respect to time, which results in
\begin{equation*} \label{est_alpha}
\begin{aligned}
\int_0^t (\alpha_h \ddot{\psi}_h, \dot{\tilde{\psi}}_h)_{L^2}\,  \textup{d}s  \, =& \, \tfrac{b}{c^2}\int_0^t\|\sqrt{\alpha_h}\, \ddot{\psi}_h\|^2_{L^2}\, \textup{d}s+\int_0^t(\alpha_h \ddot{\psi}_h, \dot{\psi}_h)_{L^2}\, \textup{d}s\\
=& \, \int_0^t \left(\tfrac{b}{c^2}\|\sqrt{\alpha_h}\, \ddot{\psi}_h\|^2_{L^2}- \tfrac12 \|\sqrt{\dot{\alpha}_{h}} \dot{\psi}_h \|^2_{L^2}\right)\, \textup{d}s \\
&+\tfrac12 \|\sqrt{\alpha_h(s)}\dot{\psi}_h(s)\|^2_{L^2} \Bigl\vert_0^t.
\end{aligned}
\end{equation*}
From here, we can further estimate the last term to obtain
\begin{equation*} 
\begin{aligned}
\int_0^t (\alpha_h \ddot{\psi}_h, \dot{\tilde{\psi}}_h)_{L^2}\,  \textup{d}s 
\geq& \, \tfrac{b}{c^2}\int_0^t\|\sqrt{\alpha_h}\, \ddot{\psi}_h\|^2_{L^2}\, \textup{d}s+\tfrac12\|\sqrt{\alpha_h(s)}\dot{\psi}_h(s)\|^2_{L^2} \Bigl\vert_0^t \\
&- \tfrac12\int_0^t \max_{x \in \Omega} \left|\frac{\dot{\alpha}_{h}(x)}{\alpha_h(x)} \right| \cdot \|\sqrt{\alpha_h}\dot{\psi}_h\|^2_{L^2}\,  \textup{d}s\\
\geq&\,  \tfrac{b}{c^2}\int_0^t\|\sqrt{\alpha_h}\, \ddot{\psi}_h\|^2_{L^2}\, \textup{d}s+\tfrac12\|\sqrt{\alpha_h(s)}\dot{\psi}_h(s)\|^2_{L^2} \Bigl\vert_0^t \\
&- \tfrac{1}{2}\max_{s \in [0,T_h]}\|\sqrt{\alpha_h(s)}\dot{\psi}_h(s)\|^2_{L^2} \, \int_0^T \max_{x \in \Omega} \left|\frac{\dot{\alpha}_{h}(x)}{\alpha_h(x)} \right| \textup{d}s,
\end{aligned}
\end{equation*}
where we have additionally employed the fact that $t\leq T_h \leq T$ in the last step. By combining our previously derived estimates, we arrive at 
\begin{equation} \label{WellpLin:est_00}
\begin{aligned}
&\begin{multlined}[t] \tfrac12 \|\sqrt{\alpha_h(t)}\dot{\psi}_h(t)\|^2_{L^2}+ C_1 c^2 \|\nabla_h \tilde{\psi}_h(t)\|^2_{L^2}\\+ \tfrac{b}{c^2}\int_0^t \|\sqrt{\alpha_h}\,\ddot{\psi}_h\|^2_{L^2}\, \textup{d}s
+C_1\|\sqrt{\chi}\, \llbracket \tilde{\psi}_h(t) \rrbracket\|^2_{\mathcal{F}_h}\end{multlined} \\
\leq & \, \begin{multlined}[t] \tfrac12 \|\sqrt{\alpha_h(0)}\dot{\psi}_h(0)\|^2_{L^2}+ C_2\, \left (c^2 \|\nabla_h \tilde{\psi}_h(0)\|^2_{L^2}+\|\sqrt{\chi}\, \llbracket \tilde{\psi}_h(0) \rrbracket \,\|^2_{\mathcal{F}_h} \right)\\
 +\tfrac{1}{2}\max_{s \in [0,T_h]}\|\sqrt{\alpha}_h(s)\dot{\psi}_h(s)\|^2_{L^2} \, \int_0^T \max_{x \in \Omega} \left|\frac{\dot{\alpha}_{h}(x)}{\alpha_h(x)} \right|  \textup{d}s  \end{multlined}
\end{aligned}
\end{equation}
for all $t \in [0,T_h]$. Taking the maximum of the above estimate over $[0,T_h]$ then yields 
\begin{equation} \label{WellpLin:est_0}
\begin{aligned}
&\begin{multlined}[t] (1-\gamma)\max_{t \in [0,T_h]}\|\sqrt{\alpha_h(t)}\dot{\psi}_h(t)\|^2_{L^2}+ c^2 \max_{t \in [0,T_h]} \|\nabla_h \tilde{\psi}_h(t)\|^2_{L^2}\\ +\tfrac{b}{c^2}\int_0^{T_h}\|\sqrt{\alpha_h}\,\ddot{\psi}_h\|^2_{L^2}\, \textup{d}s
+ \max_{t \in [0,T_h]}\|\sqrt{\chi}\, \llbracket \tilde{\psi}_h(t) \rrbracket \, \|^2_{\mathcal{F}_h}\end{multlined} \\
\lesssim & \,\begin{multlined}[t] \|\sqrt{\alpha_h(0)}\dot{\psi}_h(0)\|^2_{L^2}+ c^2 \|\nabla_h \tilde{\psi}_h(0)\|^2_{L^2}+\|\sqrt{\chi}\, \llbracket \tilde{\psi}_h(0) \rrbracket \, \|^2_{\mathcal{F}_h} . \end{multlined}
\end{aligned}
\end{equation}
Since the right-hand side of \eqref{WellpLin:est_0} does not depend on $T_h$, we are allowed to extend the existence interval to $[0,T]$; i.e., we can set $T_h=T$. Uniqueness follows by linearity of the problem and the derived stability bound.
\end{proof}
~\\
\indent Before moving to the error analysis, let us discuss how to obtain a bound on $c^2|\nabla_h \psi_h(t)|_{L^2}$ and $\tfrac{b}{c^2}|\nabla_h \dot{\psi}_h(t)|_{L^2}$ from the available bound on $c^2 |\nabla_h \tilde{\psi}_h(t)|_{L^2}$.  From the energy estimate \eqref{WellpLin:final_est}, by recalling that $\tilde{\psi}_h=\psi_h+\tfrac{b}{c^2}\dot{\psi}_h$, we have
\begin{equation} \label{est_derived}
\begin{aligned}
c^2\|\nabla_h \tilde{\psi}_h(t)\|^2_{L^2}=&\, c^2\|\nabla_h \psi_h(t)\|^2_{L^2}+\tfrac{b^2}{c^2}\|\nabla_h \dot{\psi}_h(t)\|^2_{L^2}-2 c^2(\nabla_h \psi_h(t), \tfrac{b}{c^2}\nabla_h \dot{\psi}_h(t))_{L^2} \\
\lesssim& \, E[\psi_h](0)
\end{aligned}
\end{equation}
for all $t \in [0,T]$. We can rely on the Fundamental theorem of calculus to show that
\begin{equation} \label{fund_thm_est} 
\begin{aligned}
\|\nabla_h \psi_h(t)\|_{L^2}=& \left \|\int_0^t \nabla_h \dot{\psi}_{h} \, \textup{d}s +\nabla_h \psi_h(0) \right \|_{L^2} \\
\leq& \int_0^t \|\nabla_h \dot{\psi}_{h}\|_{L^2}\, \textup{d}s +\|\nabla_h \psi_h(0)\|_{L^2} \\
\leq&\, \sqrt{t} \|\nabla_h \dot{\psi}_{h}\|_{L^2(0,t; L^2)} +\|\nabla_h \psi_h(0)\|_{L^2}.
\end{aligned}
\end{equation}
By employing Young's inequality together with inequality \eqref{fund_thm_est} in estimate \eqref{est_derived}, we arrive at
\begin{equation*} \label{est_derived_1}
\begin{aligned}
& c^2\|\nabla_h \psi_h(t)\|^2_{L^2}+\tfrac{b^2}{c^2}\|\nabla_h \dot{\psi}_h(t)\|^2_{L^2}\\
\lesssim&\, 2 c^2 \|\nabla_h \psi_h(t)\|^2_{L^2}+ \tfrac{1}{2}\tfrac{b^2}{c^2}\|\dot{\psi}_h(t)\|^2_{L^2} + E[\psi_h](0)\\
\lesssim&\,  4 c^2 T \int_0^t \|\nabla_h \dot{\psi}_h\|^2_{L^2}\, \textup{d}s+4c^2\|\nabla_h \psi_h(0)\|^2_{L^2}+ \tfrac{1}{2}\tfrac{b^2}{c^2}\|\dot{\psi}_h(t)\|^2_{L^2}+ E[\psi_h](0).
\end{aligned}
\end{equation*}
Therefore, we can conclude that
\begin{equation*} \label{est_derived_2}
\begin{aligned}
 c^2\|\nabla_h \psi_h(t)\|^2_{L^2}+\tfrac{1}{2}\tfrac{b^2}{c^2}\|\nabla_h \dot{\psi}_h(t)\|^2_{L^2}
\lesssim&\, \begin{multlined}[t]  c^2 T \int_0^t \|\nabla_h \dot{\psi}_h\|^2_{L^2}\, \textup{d}s + c^2\|\nabla_h \psi_h(0)\|^2_{L^2}\\+ E[\psi_h](0). \end{multlined}
\end{aligned}
\end{equation*}
By making use of Gronwall's inequality, we obtain
\begin{equation*} \label{Gronwall_application}
\begin{aligned}
c^2\|\nabla_h \psi_h(t)\|^2_{L^2}+\tfrac{b^2}{c^2}\|\nabla_h \dot{\psi}_h(t)\|^2_{L^2}
\lesssim\,  c^2\|\nabla_h \psi_{0,h}\|^2_{L^2}+E[\psi_h](0),
\end{aligned}
\end{equation*}
where, compared to estimate \eqref{WellpLin:final_est}, now the hidden constant depends on the final time $T$.
\section{Error analysis of the linearization} \label{Sec:Semi-discrete_Linear_Error}
In this section, we derive an \emph{a priori} error estimate for the semi-discrete problem \eqref{Westervelt_ibvp_discr}. We note that we also have to take the error of the variable coefficient $\alpha$ into account to be able to later employ a fixed-point argument and prove a convergence result for the Westervelt equation as well.  
\subsection{Error estimate in the energy norm}
We decompose the approximation error by involving the interpolant as follows:
\begin{align*}
e=\underbrace{(\psi-\psi_I)}_{e_I}-\underbrace{(\psi_h-\psi_I)}_{e_h},
\end{align*}
where $\psi$ solves \eqref{Westervelt_ibvp_cont_linearized}, $\psi_h$ solves \eqref{Westervelt_lin_ibvp_discr}, and $\psi_I$ is the interpolant introduced in~\Cref{Lemma:Interpolant}. To simplify the exposition, we introduce the following auxiliary variables:
\begin{equation*}
\begin{aligned}
\tilde{e}=\tilde{\psi}-\tilde{\psi}_h, \quad \tilde{\psi}_I=\psi_I+\tfrac{b}{c^2}\dot{\psi}_I, \quad \tilde{e}_h=\tilde{\psi}_h-\tilde{\psi}_I, \quad \tilde{e}_I=\tilde{\psi}-\tilde{\psi}_I.
\end{aligned}
\end{equation*} 
In order to formulate the convergence result, we also define the energy norm
\begin{equation*}
\begin{aligned}
||\psi_h ||_{E}:=&\, \begin{multlined}[t] \left(\max_{t \in [0,T]}\|\dot{\psi}_h(t)\|^2_{L^2}+\tfrac{b}{c^2}\int_0^T \|\ddot{\psi}_h\|_{L^2}^2 \, \textup{d}s \right. \\ \left.
+ c^2\max_{t \in [0,T]} \|\nabla_h \tilde{\psi}_h(t)\|^2_{L^2}+\max_{t \in [0,T]}|\sqrt{\chi}\, \llbracket \tilde{\psi}_h(t) \rrbracket|^2_{\mathcal{F}_h}\right)^{1/2}. \end{multlined}
\end{aligned}
\end{equation*}
Thanks to Lemma~\ref{Lemma:Interpolant}, we can estimate the interpolation error in this norm by
\begin{equation}\label{error_interpolant_energy}
\begin{aligned}
\|e_I\|^2_E
\lesssim& \,\begin{multlined}[t]\, \max_{t \in [0,T]}  \sum_{\kappa \in \mathcal{T}_h} h_\kappa^{2 \mu-2} \left(c^2|\psi(t)|^2_{H^n(\kappa)}+(1+\tfrac{b^2}{c^2})|\dot{\psi}(t)|^2_{H^n(\kappa)} \right)\\
+\tfrac{b}{c^2} \int_0^T   \sum_{\kappa \in \mathcal{T}_h} h_\kappa^{2 \mu-2}| \ddot{\psi}|^2_{H^n(\kappa)}\, \textup{d}s. \end{multlined}
\end{aligned}
\end{equation}
We are now ready to state the convergence result.
\begin{theorem}\label{Thm:LinError} Let the assumptions of Theorem~\ref{Thm:LinStability} hold. Let $\psi \in C^2([0,T]; H_0^1(\Omega) \cap H^n(\Omega))$, where $n \geq 2$, be the solution of the linear initial-boundary value problem \eqref{Westervelt_ibvp_cont_linearized}. Let $\psi_h$ be the solution of the corresponding semi-discrete problem \eqref{Westervelt_lin_ibvp_discr} with the approximate initial data given by $$(\psi_{h}(0), \dot{\psi}_h(0))=((\psi_0)_I, (\psi_1)_I),$$ and the parameter $\beta$ in \eqref{stabilization} chosen sufficiently large according to Theorem~\ref{Thm:LinStability}. Then the following bound holds for the discretization error:
\begin{equation} \label{ErrorLin}
\begin{aligned}
\|\psi-\psi_h\|^2_E
\leq& \,C_{\textup{Th}\ref{Thm:LinError}} \begin{multlined}[t]\,\left\{h^{2 \mu-2} \max_{t \in [0,T]}  \sum_{\kappa \in \mathcal{T}_h} \left( |\dot{\psi}(t)|^2_{H^n(\kappa)}+|\psi(t)|^2_{H^n(\kappa)} \right) \right.\\
+h^{2 \mu-2} \int_0^T   \sum_{\kappa \in \mathcal{T}_h} \left (| \ddot{\psi}|^2_{H^n(\kappa)}+|\dot{\psi}|^2_{H^n(\kappa)}\right)\, \textup{d}s\\ \left.
+\int_0^T \|(\alpha-\alpha_h) \ddot{\psi}\|^2_{L^2} \, \textup{d}s \vphantom{ \frac{h_\kappa^{2 \mu-2}}{p^{2s-3}}} \right\}, \end{multlined}
\end{aligned}
\end{equation}
\noindent where $\mu =\min \{n, p+1\}$ and $h=\max_{\kappa \in \mathcal{T}_h} h_\kappa$, provided that $\gamma \in (0,1)$ in \eqref{alpha_t_small} is sufficiently small. The constant $C_{\textup{Th}\ref{Thm:LinError}}>0$ depends on the polynomial degree, but not  on the mesh size.
\end{theorem}
\begin{proof}
We begin the proof by observing that $\psi$ satisfies the weak form \eqref{Westervelt_lin_ibvp_discr} when $\alpha_h=\alpha$. Therefore, we can see the error $e=\psi-\psi_h$ as the solution of the following problem: 
\begin{equation}
\begin{aligned} \label{ibvp_error}
\begin{multlined}[t]\displaystyle 
(\alpha_h\ddot{e}, v_h)_{L^2}+a_h(\tilde{e}, v_h) =- ((\alpha-\alpha_h)\ddot{\psi}, v_h)_{L^2}  \end{multlined}
\end{aligned}
\end{equation}
for all  $v_h \in V_h$ and all time  $t\in (0,T]$, with $$(e(0), \dot{e}(0))=(\psi_0-\psi_{0, h},\, \psi_1-\psi_{1, h}).$$ By involving the interpolant, we can then rewrite equation \eqref{ibvp_error} as 
\begin{equation} \label{ibvp_error_rewritten}
\begin{aligned}
 (\alpha_h\ddot{e}_h, v_h)_{L^2}+a_h(\tilde{e}_h, v_h) 
= (\alpha_h\ddot{e}_I+(\alpha-\alpha_h)\ddot{\psi}, v_h)_{L^2}+a_h(\tilde{e}_I, v_h) 
\end{aligned}
\end{equation}
for $t \in (0,T]$ and all $v_h \in V_h$. We next test equation \eqref{ibvp_error_rewritten} with $v_h=\dot{\tilde{e}}_h \in V_h$ and estimate the resulting terms. By treating all the terms arising from the left-hand side of \eqref{ibvp_error_rewritten} as in the proof of Theorem~\ref{Thm:LinStability},  we arrive at the following counterpart of estimate \eqref{WellpLin:est_00} for the energy of $e_h$ at time $t$:
\begin{equation} \label{ErrorLin:1st_est}
\begin{aligned}
&\,\begin{multlined}[t] \|\sqrt{\alpha_h(t)}\, \dot{e}_h(t)\|^2_{L^2}+\tfrac{b}{c^2}\int_0^t \|\sqrt{\alpha_h}\, \ddot{e}_h\|_{L^2}^2 \, \textup{d}s  + c^2\|\nabla_h \tilde{e}_h(t)\|^2_{L^2}
	+\|\sqrt{\chi}\, \llbracket \tilde{e}_h(t) \rrbracket \,\|^2_{\mathcal{F}_h} \end{multlined}\\
\lesssim& \,\begin{multlined}[t] \left|\int_0^t \left\{ (\alpha_h\ddot{e}_I+(\alpha-\alpha_h)\ddot{\psi},\, \dot{\tilde{e}}_h)_{L^2}+a_h(\tilde{e}_I, \dot{\tilde{e}}_h)\right\} \, \textup{d}s \right|
+\gamma \max_{s \in [0,t]}\|\sqrt{\alpha_h(s)}\dot{e}_h(s)\|^2_{L^2},
 \end{multlined}
\end{aligned}
\end{equation}
where we have additionally used that $e_h(0)=\dot{e}_h(0)=0$ due to our choice of the approximate initial data. By employing H\"older's and then Young's inequality with $\varepsilon_1>0$ , we obtain
\begin{equation*}
\begin{aligned}
\qquad &\int_0^t (\alpha_h\, \ddot{e}_I+(\alpha-\alpha_h)\ddot{\psi},\, \dot{e}_h+\tfrac{b}{c^2}\ddot{e}_h)_{L^2}\, \textup{d}s\\
\leq& \, \int_0^t  (\|\sqrt{\alpha_h}\ddot{e}_I\|_{L^2}+ \left \|\tfrac{\alpha-\alpha_h}{\sqrt{ \alpha_h}}\ddot{\psi}\right \|_{L^2})(\|\sqrt{\alpha_h}\dot{e}_h\|_{L^2}+\tfrac{b}{c^2}\|\sqrt{\alpha_h}\ddot{e}_h \|_{L^2})\, \textup{d}s \\
\leq&\,\begin{multlined}[t] \tfrac{1}{4\varepsilon_1}(1+\tfrac{b}{c^2})\int_0^t \left(\alpha_1 \|\ddot{e}_{I}\|^2_{L^2}+\alpha_0^{-1}\|(\alpha-\alpha_h)\ddot{\psi}\|^2_{L^2} \right)\, \textup{d}s\\+ 2\varepsilon_1 \int_0^t \|\sqrt{\alpha_h}\dot{e}_h\|^2_{L^2}\, \textup{d}s 
+2\varepsilon_1 \tfrac{b}{c^2} \int_0^t \|\sqrt{\alpha_h}\ddot{e}_h\|^2_{L^2}\, \textup{d}s \end{multlined}
\end{aligned}
\end{equation*}
for all $t \in [0,T]$. We note that the following useful estimate holds for the bilinear form $a_h(\cdot, \cdot)$:
\begin{equation} \label{ineq_Ah}
\begin{aligned}
|a_h(\phi, v_h)| \lesssim&\, \begin{multlined}[t] c^2\|\nabla_h \phi\|_{L^2}\|\nabla_h v \|_{L^2}+\|\chi^{-1/2}\llbrace c^2 \nabla_h \phi \rrbrace \,\|_{\mathcal{F}_h} \|\sqrt{\chi} \llbracket v_h \rrbracket \, \|_{\mathcal{F}_h}\\
+c\tfrac{1}{\sqrt{\beta}}\| \nabla_h v_h\|_{L^2} \|\sqrt{\chi} \llbracket \phi \rrbracket \, \|_{\mathcal{F}_h} +\|\sqrt{\chi} \llbracket \phi \rrbracket \,\|_{\mathcal{F}_h}\|\sqrt{\chi} \llbracket v_h \rrbracket \,\|_{\mathcal{F}_h}, \end{multlined}
\end{aligned}
\end{equation}
which we will use below with the choice of $\phi \in \{\tilde{e}_I, \dot{\tilde{e}}_I\}$. The term $\|\chi^{-1/2}\llbrace c^2 \nabla_h \phi \rrbrace\|_{\mathcal{F}_h}$ appearing in \eqref{ineq_Ah} is the reason why we need the bound \eqref{interpolant_error_2} on the interpolant error. To estimate the $a_h(\tilde{e}_I, \dot{\tilde{e}}_h)$ term in \eqref{ErrorLin:1st_est}, we employ integration by parts with respect to time and then twice inequality \eqref{ineq_Ah},
\begin{equation*}
\begin{aligned}
\int_0^t a_h(\tilde{e}_I, \dot{\tilde{e}}_h) \, \textup{d}s=&\, a_h(\tilde{e}_I(t), \tilde{e}_h(t)) -\int_0^t a_h(\dot{\tilde{e}}_I, \tilde{e}_h) \, \textup{d}s \\
\lesssim &\, \begin{multlined}[t] c^2\|\nabla_h \tilde{e}_I(t)\|_{L^2}|\nabla_h \tilde{e}_h(t)|_{L^2}+\|\chi^{-1/2}\llbrace c^2 \nabla_h \tilde{e}_I(t) \rrbrace \|_{\mathcal{F}_h} \|\sqrt{\chi} \llbracket \tilde{e}_h(t) \rrbracket \|_{\mathcal{F}_h}\\[1mm]
+c\tfrac{1}{\sqrt{\beta}}\| \nabla_h \tilde{e}_h(t)\|_{L^2} \|\sqrt{\chi} \llbracket \tilde{e}_I(t) \rrbracket \|_{\mathcal{F}_h}+\|\chi^{1/2}\llbracket \tilde{e}_I(t) \rrbracket\|_{\mathcal{F}_h} \|\chi^{1/2}\llbracket \tilde{e}_h(t) \rrbracket\|_{\mathcal{F}_h} \\
+\int_0^t \left(c^2\|\nabla_h \dot{\tilde{e}}_I\|_{L^2}\|\nabla_h \tilde{e}_h\|_{L^2}+\|\chi^{-1/2}\llbrace c^2 \nabla_h \dot{\tilde{e}}_I \rrbrace\|_{\mathcal{F}_h} \|\sqrt{\chi} \llbracket \tilde{e}_h \rrbracket \|_{\mathcal{F}_h} \right. \\ \left.
+c\tfrac{1}{\sqrt{\beta}}\| \nabla_h \tilde{e}_h\|_{L^2} \|\sqrt{\chi} \llbracket \dot{\tilde{e}}_I \rrbracket \|_{\mathcal{F}_h} +\|\sqrt{\chi}\llbracket \dot{\tilde{e}}_I \rrbracket\|_{\mathcal{F}_h} \|\sqrt{\chi}\llbracket \tilde{e}_h \rrbracket\|_{\mathcal{F}_h} \right)\, \textup{d}s.
\end{multlined} 
\end{aligned}
\end{equation*}
From here by Young's inequality with $\varepsilon_2 \in (0, \varepsilon_1)$, we have 
\begin{equation*}
\begin{aligned}
\int_0^t a_h(\tilde{e}_I, \dot{\tilde{e}}_h) \, \textup{d}s
\lesssim&\, \begin{multlined}[t] \varepsilon_1 \left(c^2(1+\tfrac{1}{\beta})\|\nabla_h \tilde{e}_h(t)\|^2_{L^2}
+ \|\sqrt{\chi} \llbracket \tilde{e}_h(t) \rrbracket \|^2_{\mathcal{F}_h}
\right) \\
+\varepsilon_2 \int_0^t \left(c^2(1+\tfrac{1}{\beta})\|\nabla \tilde{e}_h\|_{L^2}^2+\|\sqrt{\chi} \llbracket \tilde{e}_h \rrbracket \|^2_{\mathcal{F}_h}\right)\, \textup{d}s+ \tfrac{1}{4\varepsilon_2} \, \bar{E}[e_I](t),
\end{multlined} 
\end{aligned}
\end{equation*}
where the modified energy of the interpolant error is given by
	\begin{equation*} \label{energy_interpolant}
	\begin{aligned}
			\bar{E}[e_I](t):=&\, \begin{multlined}[t]  c^2 \|\nabla_h \tilde{e}_I(t)\|^2_{L^2}+\|\sqrt{\chi} \llbracket \tilde{e}_I(t) \rrbracket \|^2_{\mathcal{F}_h}
	+\|\chi^{-1/2}\llbrace c^2 \nabla_h \tilde{e}_I(t) \rrbrace\|^2_{\mathcal{F}_h}
	 \\
	+\int_0^t \left(c^2\|\nabla_h \dot{\tilde{e}}_I\|^2_{L^2}+\|\sqrt{\chi} \llbracket \dot{\tilde{e}}_I \rrbracket \|^2_{\mathcal{F}_h}
+ \|\chi^{-1/2}\llbrace c^2 \nabla_h \dot{\tilde{e}}_I \rrbrace\|^2_{\mathcal{F}_h} \right)\, \textup{d}s. 	\end{multlined}
	\end{aligned}
	\end{equation*} 	
 We fix $\varepsilon_1>0$ sufficiently small and include the derived bounds in estimate \eqref{ErrorLin:1st_est}, from which we immediately have
		\begin{equation} \label{ErrorLin:energy_est1}
		\begin{aligned}
	E[e_h](t) 
		\lesssim& \,\begin{multlined}[t] \varepsilon_2 T \max_{t \in [0,T]}\left(\|\sqrt{\alpha_h(t)}\dot{e}_h(t)\|^2_{L^2}+ c^2(1+\tfrac{1}{\beta}) \|\nabla_h \tilde{e}_h(t)\|^2_{L^2}\right. \\ \left.+\|\sqrt{\chi} \llbracket \tilde{e}_h(t) \rrbracket \|^2_{\mathcal{F}_h}\right) 
		+(1+\tfrac{b}{c^2})\, \alpha_0^{-1} \int_0^t \|(\alpha-\alpha_h)\ddot{\psi}\|^2_{L^2} \, \textup{d}s \\+\tfrac{1}{4 \varepsilon_2}\left(\bar{E}[e_I](t)+(1+\tfrac{b}{c^2})\alpha_1\int_0^t  \|\ddot{e}_{I}\|^2_{L^2}\, \textup{d}s \right)\\
		+\gamma \max_{t \in [0,T]}\|\sqrt{\alpha_h(t)}\dot{e}_h(t)\|^2_{L^2}
		\end{multlined}
		\end{aligned}
		\end{equation}
	for all $t \in [0,T]$. Above, we have also employed the inequality \[\int_0^t \|v\|^2_{L^2} \, \textup{d}s \leq T \max_{t \in [0,T]} \|v(t)\|^2_{L^2},\] which holds for functions $v \in C([0,T]; L^2(\Omega))$.
	By possibly decreasing $\varepsilon_2$ and $\gamma$, and then taking the maximum over $[0,T]$ of \eqref{ErrorLin:energy_est1}, we obtain
	\begin{equation} \label{ErrorLin:energy_est3}
	\begin{aligned}
\max_{t \in [0,T]}	E[e_h](t) 
	\lesssim& \,\begin{multlined}[t] \, (1+\tfrac{b}{c^2})\, \alpha_0^{-1} \int_0^T \|(\alpha-\alpha_h)\ddot{\psi}\|^2_{L^2} \, \textup{d}s+\max_{t \in [0,T]}\bar{E}[e_I](t)\\ +(1+\tfrac{b}{c^2})\alpha_1\int_0^T \|\ddot{e}_{I}\|^2_{L^2}\, \textup{d}s.
	\end{multlined}
	\end{aligned}
	\end{equation}
Recalling also the properties of the interpolant stated in~\Cref{Lemma:Interpolant} leads to
\begin{equation*}
\begin{aligned}
\|e_h\|^2_{E} 
\lesssim& \,\begin{multlined}[t]\,h^{2 \mu-2}  \max_{t \in [0,T]}  \sum_{\kappa \in \mathcal{T}_h}\left( c^2|\psi(t)|^2_{H^n(\kappa)}+\tfrac{b^2}{c^2}|\dot{\psi}(t)|^2_{H^n(\kappa)}\right)\\
+h^{2 \mu-2} \int_0^T   \sum_{\kappa \in \mathcal{T}_h} \left ((\alpha_1+\tfrac{b}{c^2}\alpha_1+\tfrac{b^2}{c^2})| \ddot{\psi}|^2_{H^n(\kappa)}+c^2|\dot{\psi}|^2_{H^n(\kappa)}\right)\, \textup{d}s\\
 +(1+\tfrac{b}{c^2})\, \alpha_0^{-1} \int_0^T \|(\alpha-\alpha_h)\ddot{\psi}\|^2_{L^2} \, \textup{d}s. \end{multlined}
\end{aligned}
\end{equation*}
Together with estimate \eqref{error_interpolant_energy} for $\|e_I\|_E$, this yields the desired bound \eqref{ErrorLin} for the discretization error. We note that the constant $C_{\textup{Thm}\ref{Thm:LinError}}$ in the final estimate has the form
\begin{align*}
C_{\textup{Thm}\ref{Thm:LinError}}=O\left(c^2+\tfrac{b^2}{c^2}+\tfrac{b}{c^2}+1\right),
\end{align*}
and so it does not degenerate as $b$ approaches zero.
\end{proof}

\section{Analysis of the nonlinear model} \label{Sec:Semi-discrete_Nonlinear_Error}
Our next aim is to analyze the semi-discretization of the Westervelt equation given by \eqref{Westervelt_ibvp_discr}. To this end, we will rely on our analysis of the linearized problem together with a fixed-point argument.
\begin{theorem} \label{Thm:ErrorNonlinear}
Let $c>0$, $b \geq 0$, and $k \in \mathbb{R}$. Let the final time $T'\leq T$ be such that the initial-boundary value problem \eqref{Westervelt_ibvp_cont} for the Westervelt equation has a solution $$\psi \in C^2([0,T']; H_0^1(\Omega) \cap H^n(\Omega)),\ \text{where } n >1+d/2, $$ for which it holds that
\begin{equation*}
\begin{aligned}
0< \alpha_0 \leq 1-2k\dot{\psi} \leq \alpha_1  \qquad \text{in } \overline{\Omega} \times [0,T']
\end{aligned}
\end{equation*}
for some $\alpha_0$, $\alpha_1>0$.  Assume that the polynomial degree $p \geq 2$ and that the approximate initial conditions are given by 
\begin{equation} \label{approx_initial_conditions}
\begin{aligned}
(\psi_{h}(0), \dot{\psi}_h(0))=((\psi_0)_I, (\psi_1)_I).
\end{aligned}
\end{equation}
Then for sufficiently small $h=\displaystyle \max_{\kappa \in \mathcal{T}_h} h_\kappa$ and
\begin{align*}
M(\psi)=\begin{multlined}[t] \max_{t \in [0,T']} |\psi(t)|^2_{H^n}+\max_{t \in [0,T']}|\dot{\psi}(t)|^2_{H^n}+\max_{t \in [0,T']}\|\dot{\psi}(t)\|^2_{L^\infty}\\+\int_0^{T'}(|\ddot{\psi}|^2_{H^n}+\|\ddot{\psi}\|^2_{L^\infty}+|\dot{\psi}|^2_{H^n}) \,\textup{d}s, \end{multlined}
\end{align*}
 the corresponding semi-discrete problem \eqref{Westervelt_ibvp_discr} for the Westervelt equation has a unique solution $\psi_h \in C^2([0,T']; V_h)$ that satisfies the following error bound: 
\begin{equation*} \label{error_est_Westervelt}
\begin{aligned}
||\psi-\psi_h||^2_{E} \leq \, \begin{multlined}[t] C_{\textup{Th}\, \ref{Thm:ErrorNonlinear}}\, h^{2 \mu-2} \,\sum_{\kappa \in \mathcal{T}_h} \left\{ \max_{t \in [0,T']}   \left( |\dot{\psi}(t)|^2_{H^n(\kappa)}+|\psi(t)|^2_{H^n\kappa)} \right) \right.\\ \left.
+\int_0^{T'}  \left (| \ddot{\psi}|^2_{H^n(\kappa)}+|\dot{\psi}|^2_{H^n(\kappa)}\right)\, \textup{d}s\right\} , \end{multlined}
\end{aligned}
\end{equation*}
with $\mu =\min \{n, p+1\}$, provided that the parameter $\beta$ in \eqref{stabilization} is sufficiently large. The  constant $C_{\textup{Th}\ref{Thm:ErrorNonlinear}}>0$ depends on $M(\psi)$ and on the polynomial degree, but not on the mesh size.
\end{theorem}
\begin{proof}
We conduct the proof by employing the Banach fixed-point theorem. Therefore, we first need to define a fixed-point map. We begin by introducing the set
\begin{equation*} \label{set_Bh}
\begin{aligned}
\mathcal{B}_h=\left \{w_h \in C^2([0,T']; V_h):\vphantom{||\psi-w_h |||^2_{\mathcal{E}}} \right. \ & \left. ||\psi-w_h ||^2_{E} \leq \,  C_{\textup{Th}\ref{Thm:ErrorNonlinear}}\, h^{2 \mu-2}  \sum_{\kappa \in \mathcal{T}_h}  \max_{t \in [0,T']}  \mathcal{E}_\kappa[\psi](t) \, \right. \\
& \left. (w_h(0), \dot{w}_{h}(0))=(\psi_{0, h}, \psi_{1, h})\, \vphantom{|||\psi-w_h |||^2_{\mathcal{E}}} \right\},
\end{aligned}
\end{equation*}
where we have used the notation
\begin{equation*} \label{E_kappa}
\begin{aligned}
\mathcal{E}_\kappa[\psi](t)=& \, |\dot{\psi}(t)|^2_{H^n(\kappa)}+|\psi(t)|^2_{H^n(\kappa)}+\int_0^t \left (| \ddot{\psi}|^2_{H^n(\kappa)}+|\dot{\psi}|^2_{H^n(\kappa)}\right)\, \textup{d}s.
\end{aligned}
\end{equation*}
The constant $C_{\textup{Th}\ref{Thm:ErrorNonlinear}}>0$ will be made precise below. \\[3mm]
\noindent {\bf Step 1: Defining the fixed-point map.}\\[2mm]
\noindent  For $w_h \in \mathcal{B}_h$, we then define the operator $\mathcal{J}: w_h \mapsto \psi_h$, where $\psi_h$ solves
\begin{align} \label{Westervelt_lin_ibvp_discr_fixedpoint}
\begin{cases}
\displaystyle ((1-2k \dot{w}_h)\ddot{\psi}_h,\, v_h)_{L^2}+a_h(\tilde{\psi}_h, v_h) = 0   \qquad
\text{for all } v_h \in V_h, \ t\in (0,T'), \\[3mm]
(\psi_h(0), \dot{\psi}_{h}(0))=(\psi_{0, h}, \psi_{1, h}).
\end{cases}
\end{align}
The operator $\mathcal{J}$ is well-defined thanks to the well-posedness result of Theorem~\ref{Thm:LinStability}, whose assumptions we verify below. We note that the set $\mathcal{B}_h$ is non-empty because the global interpolant $\psi_I$ is in $\mathcal{B}_h$. Moreover, $\mathcal{B}_h$ is closed with respect to topology induced by $\|\cdot\|_{E}$.\\[3mm]
\noindent \textbf{Step 2: The self-mapping property.}\\[2mm]
\noindent We next want to verify that $\mathcal{J}(\mathcal{B}_h) \subset \mathcal{B}_h$. Let $w_h \in \mathcal{B}_h$. To show that $\psi_h=\mathcal{J}(w_h) \in B_h$, we rely on Theorems~\ref{Thm:LinStability} and ~\ref{Thm:LinError}, which guarantee stability and convergence for the linearized problem. We choose the variable coefficient in the linear model to be $\alpha_h=1-2k \dot{w}_h$ and check that all the assumptions of these theorems are satisfied. In particular, we have to justify the non-degeneracy assumption \eqref{assumptions_alpha_h} and the smallness of $\dot{\alpha}_h/\alpha_h$ in \eqref{alpha_t_small}. \\
\indent  We already know that $\alpha_h=1-2k \dot{w}_h \in C^1([0,T']; V_h)$. We next rely on the inverse estimate given in~\Cref{lem:inversepoly} and properties of the interpolant stated in~\Cref{Lemma:Interpolant} to verify the non-degeneracy assumption \eqref{assumptions_alpha_h} on $\alpha_h$. \\[3mm]
\indent \emph{The coefficient $\alpha_h$ does not degenerate.} Fix $t \in [0,T']$. We can pick an element $\hat{\kappa} \in \mathcal{T}_h$, such that
\begin{equation} \label{nondegeneracy_bound_0}
\begin{aligned}
\max_{x \in \overline{\Omega}}|\dot{w}_h(x, t)|= \,\max_{x \in \hat{\kappa}}|\dot{w}_h(x, t)|.
\end{aligned}
\end{equation}
By involving the local interpolant and then relying on the inverse estimate, we find that
\begin{equation} \label{nondegeneracy_bound_1}
\begin{aligned}
\max_{x \in \overline{\Omega}}|\dot{w}_h(x, t)|^2\ =\qquad & \,\max_{x \in \hat{\kappa}}  |\dot{w}_h(x, t)|^2\\
\stackrel{\text{triangle ineq.}}{\lesssim}&   \|\dot{\psi}_{I, \hat{\kappa}}(t)-\dot{w}_h(t)\|^2_{L^\infty(\hat{\kappa})}+\|\dot{\psi}_{I, \hat{\kappa}}(t)\|^2_{L^\infty(\hat{\kappa})}  \\ 
\stackrel{\text{inverse est.} }{\lesssim}& \,  h_{\hat{\kappa}}^{-d} \, \|\dot{\psi}_{I, \hat{\kappa}}(t)-\dot{w}_h(t)\|^2_{L^2(\hat{\kappa})}+\|\dot{\psi}_{I, \hat{\kappa}}(t)\|^2_{L^\infty(\hat{\kappa})}  \\
\stackrel{\text{triangle ineq.}}{\lesssim}&\, \begin{multlined}[{t}]  h_{\hat{\kappa}}^{-d} \left(\|\dot{\psi}_{I, \hat{\kappa}}(t)-\dot{\psi}(t)\|^2_{L^2(\hat{\kappa})}+\|\dot{\psi}(t)-\dot{w}_h(t)\|^2_{L^2(\hat{\kappa})}\right)\\
+\|\dot{\psi}_{I, \hat{\kappa}}(t)\|^2_{L^\infty(\hat{\kappa})}.\end{multlined}   
\end{aligned}
\end{equation}
We can estimate the last three terms on the right-hand side of \eqref{nondegeneracy_bound_1} by employing the stability and approximation properties of the interpolant, and the fact that $w_h \in \mathcal{B}_h$. By doing so, we obtain
\begin{equation*} 
\begin{aligned}
\max_{x \in \overline{\Omega}} |\dot{w}_h(x, t)|^2
\lesssim&\, \begin{multlined}[{t}]   \, h_{\hat{\kappa}}^{2 \mu-d}|\dot{\psi}(t)|^2_{H^s(\hat{\kappa})}
+ h_{\hat{\kappa}}^{-d} h^{2\mu-d}\,C_{\textup{Th}\ref{Thm:ErrorNonlinear}} \sum_{\kappa \in \mathcal{T}_h}\,\max_{t \in [0,T']}  \mathcal{E}_\kappa[\psi](t)+ \|\dot{\psi}_{I, \hat{\kappa}}(t)\|^2_{L^\infty(\hat{\kappa})} \end{multlined}\\
\lesssim&\, \begin{multlined}[t] \, h^{2\mu-d}  |\dot{\psi}(t)|^2_{H^s(\Omega)}
+ C_{\textup{Th}\ref{Thm:ErrorNonlinear}} \left(\frac{h}{h_{\hat{\kappa}}}\right)^d \,  h^{2 \mu-2-d}\, \sum_{\kappa \in \mathcal{T}_h}\, \max_{t \in [0,T']} \mathcal{E}_\kappa [\psi](t) \\+ \|\dot{\psi}_{I, \hat{\kappa}}(t)\|^2_{L^\infty(\hat{\kappa})} ,\end{multlined}\\
\end{aligned}
\end{equation*}
recalling that due to our assumptions, we have $\mu=\min \{n, p+1\} > 1+d/2$. By using the $L^\infty$ stability of the interpolant and the assumption on quasi-uniformity of the mesh, we infer
\begin{equation} \label{nondegeneracy_bound_4}
\begin{aligned}
 \max_{x \in \overline{\Omega}}|\dot{w}_h(x, t)|^2
\leq\, \begin{multlined}[{t}] C_1 M(\psi) \left(h^{2\mu-d}
+ C_{\textup{Th}\ref{Thm:ErrorNonlinear}}(M(\psi))\, h^{2\mu-2-d}
+1 \right)\, :=m\end{multlined}
\end{aligned}
\end{equation}
for every $t \in [0,T']$, where the constant $C_1>0$ above does not depend on the mesh size. We refer to \eqref{constant_Thm3} below for the exact form of $C_{\textup{Th}\ref{Thm:ErrorNonlinear}}(M(\psi))$.  By taking the maximum over $t \in [0,T']$ in  \eqref{nondegeneracy_bound_4},  we further have \[
\|\dot{w}_h\|_{C(\overline{\Omega} \times [0,T'])} \leq \sqrt{m} .
\]
 We then choose $M(\psi)$ and $h$ sufficiently small so that
$$0<\alpha_0 \leq 1-2|k|\sqrt{m} \leq \alpha_h=1-2k \dot{w}_h \leq 1+2|k|\sqrt{m} \leq \alpha_1$$
in $\overline{\Omega} \times [0,T']$.\\[3mm]
\indent \emph{The quotient $\dot{\alpha}_h/\alpha_h$ is sufficiently small}. The assumption on $\dot{\alpha_h}/ \alpha_h$ in Theorems~\ref{Thm:LinStability} and~\ref{Thm:LinError} can be verified as follows:
\begin{equation*}
\begin{aligned}
\left\|\frac{\dot{\alpha}_h}{\alpha_h} \right\|_{L^1 L^\infty}=\left\|\frac{-2k \ddot{w}_h}{1-2k\dot{w}_h} \right\|_{L^1 L^\infty} \leq&\, \frac{2 |k|}{1-2|k| \sqrt{m}}\|\ddot{w}_h\|_{L^1 L^\infty} \\
\leq&\,  \frac{2 |k|}{1-2|k| \sqrt{m}}\sqrt{T'}\|\ddot{w}_h\|_{L^2 L^\infty}.
\end{aligned}
\end{equation*}
We can bound $\|\ddot{w}_h\|_{L^2 L^\infty}$ in a similar fashion as \eqref{nondegeneracy_bound_0}--\eqref{nondegeneracy_bound_4} by relying on the interpolant and inverse estimates. In particular, we have
\begin{equation*} 
\begin{aligned}
\displaystyle \int_0^{T'} \max_{x \in \overline{\Omega}}|\ddot{w}_h(s)|^2 \, \textup{d}s
\lesssim& \, \begin{multlined}[t]   \, h^{2\mu-d} \displaystyle \int_0^{T'}|\ddot{\psi}(s)|^2_{H^s(\Omega)}\, \textup{d}s  
+ C_{\textup{Th}\ref{Thm:ErrorNonlinear}} h^{2\mu-2-d} \,\sum_{\kappa \in \mathcal{T}_h}\, \max_{t \in [0,T']} 
 \mathcal{E}_\kappa [\psi](t) \\
 + \displaystyle \int_0^{T'}\|\ddot{\psi}_{I, \hat{\kappa}}(s)\|^2_{L^\infty(\hat{\kappa})}\, \textup{d}s,
\end{multlined}
\end{aligned}
\end{equation*}
from which it follows that
\begin{equation*} \label{wh_tt_bound}
\begin{aligned}
\displaystyle \int_0^{T'} \max_{x \in \overline{\Omega}}|\ddot{w}_h(s)|^2 \, \textup{d}s
\leq C_2 m,
\end{aligned}
\end{equation*}
where the constant $C_2>0$  does not depend on the mesh size. Therefore, for \Cref{Thm:LinError} to hold, we need that
\begin{equation*}
\begin{aligned}
\gamma:=\frac{2|k|}{1-2|k|\sqrt{m}}\sqrt{T'C_2 m}
\end{aligned}
\end{equation*}
is sufficiently small, which we can achieve by decreasing $M(\psi)$ and $h$. \\[3mm]
\indent \emph{Choosing the constant $C_{\textup{Th}\ref{Thm:ErrorNonlinear}}$ so that $\mathcal{F}$ is a self-mapping}. We have therefore verified all the assumptions of Theorems~\ref{Thm:LinStability} and~\ref{Thm:LinError}. On account of~\Cref{Thm:LinError} and the resulting error estimate \eqref{ErrorLin}, we conclude that problem \eqref{Westervelt_lin_ibvp_discr_fixedpoint} has a unique solution $\psi_h \in C^2([0,T']; V_h)$ that satisfies
\begin{equation*}
\begin{aligned}
 ||\psi-\psi_h ||^2_{E}
  \leq& \,\begin{multlined}[t] C_{\textup{Th} \ref{Thm:LinError}} \left(h^{2\mu-2}\,\sum_{\kappa \in \mathcal{T}_h} \max_{t \in [0,T']}  \mathcal{E}_\kappa[\psi](t)\right. \left. + \max_{t \in [0,T']}\|\alpha(t)-\alpha_h(t)\|^2_{L^2}\int_0^{T'} \|\ddot{\psi}\|^2_{L^\infty(\Omega)} \, \textup{d}s\right). \end{multlined}
\end{aligned}
\end{equation*}
Noting that the error in $\alpha$ can be estimated according to 
\begin{equation*} \label{error_alpha}
\begin{aligned}
\max_{t \in [0,T']} \|\alpha(t)-\alpha_h(t)\|^2_{L^2} =& \,\max_{t \in [0,T']} \|2k\dot{\psi}(t)-2k\dot{w}_h(t)\|^2_{L^2}\\ \leq& \, 4k^2 C_{\textup{Th}\ref{Thm:ErrorNonlinear}}\, h^{2\mu-2} \displaystyle \sum_{\kappa \in \mathcal{T}_h} \max_{t \in [0,T']}  \mathcal{E}_\kappa [\psi](t),
\end{aligned}
\end{equation*}
since $w_h \in \mathcal{B}_h$, we obtain
\begin{equation*}
\begin{aligned}
||\psi-\psi_h ||^2_{E} \leq& \,  C_{\textup{Th} \ref{Thm:LinError}}\,(1+4k^2 C_{\textup{Th}\ref{Thm:ErrorNonlinear}} M(\psi)) h^{2\mu-2} \sum_{\kappa \in \mathcal{T}_h} \max_{t \in [0,T']}  \mathcal{E}_\kappa [\psi](t).
\end{aligned}
\end{equation*}
For sufficiently small $M(\psi)$ such that $1-4k^2 C_{\textup{Th} \ref{Thm:LinError}}M(\psi) >0$, we can choose $C_{\textup{Th}\ref{Thm:ErrorNonlinear}}$ as
\begin{equation} \label{constant_Thm3}
\begin{aligned}
C_{\textup{Th}\ref{Thm:ErrorNonlinear}} :=\frac{C_{\textup{Th} \ref{Thm:LinError}}}{1-4k^2 C_{\textup{Th} \ref{Thm:LinError}} M(\psi)}.
\end{aligned}
\end{equation}
This choice of the constant $C_{\textup{Th}\ref{Thm:ErrorNonlinear}}$ implies that $$||\psi-\psi_h ||^2_{E} \leq \,  C_{\textup{Th}\ref{Thm:ErrorNonlinear}} h^{2 \mu-2} \sum_{\kappa \in \mathcal{T}_h} \max_{t \in [0,T']}  \mathcal{E}_\kappa [\psi](t);$$ in other words, $\psi_h \in \mathcal{B}_h$. \\

\noindent \textbf{Step 3: Contractivity.}\\[2mm]
\noindent To prove that the operator $\mathcal{J}$ is strictly contractive, take $w_h^{(1)}$, $w_h^{(2)} \in \mathcal{B}_h$ and set $\psi_h^{(1)}=\mathcal{J} (w_h^{(1)})$, $\psi_h^{(2)}=\mathcal{J} (w_h^{(2)}) \in \mathcal{B}_h$. Denote $W_h=w_h^{(1)}-w_h^{(2)}$. Then the difference $\Psi_h=\psi_h^{(1)}-\psi_h^{(2)}$ satisfies the problem
 \begin{align} \label{ibvp_contractivity}
 ((1-2k \dot{w}_h^{(1)})\ddot{\Psi}_h, v_h)_{L^2}+a_h(\tilde{\Psi}_h, v_h)
 = (2k\, \dot{W}_h\, \ddot{\psi}^{(2)}_h, v_h)_{L^2}, \qquad v_h \in V_h
 \end{align}
for all time, with zero initial data. This equation corresponds to equation \eqref{ibvp_error} satisfied by the approximation error in the proof of~\Cref{Thm:LinError}. Therefore, testing \eqref{ibvp_contractivity} with $\dot{\tilde{\Psi}}_h$ and proceeding analogously to the proof of~\Cref{Thm:LinError}  results in the estimate 
\begin{equation*}
\begin{aligned}
||\Psi_h||^2_{E} \lesssim&\, 4k^2 (1+\tfrac{b}{c^2})\alpha_0^{-1}\left \{ \int_0^{T'} \|\ddot{\psi}_h^{(2)}\|^2_{L^\infty(\Omega)}\, \textup{d}s \right \} \max_{t \in [0,T']}|\dot{W}_h(t)|^2_{L^2}.
\end{aligned}
\end{equation*}
This inequality then corresponds to estimate \eqref{ErrorLin:energy_est3} in the proof of~\Cref{Thm:LinError} if we formally set the interpolant error to zero. We further have
\begin{equation} \label{contractivity}
\begin{aligned}
||\Psi_h||^2_{E} \lesssim\,4k^2 (1+\tfrac{b}{c^2})\alpha_0^{-1} \left \{ \int_0^{T'} \|\ddot{\psi}_h^{(2)}\|^2_{L^\infty(\Omega)}\, \textup{d}s \right \} ||W_h||^2_{E}.
\end{aligned}
\end{equation}
We can bound $\|\ddot{\psi}_h^{(2)}\|^2_{L^2 L^\infty}$ by proceeding in the same way as in \eqref{nondegeneracy_bound_1}--\eqref{nondegeneracy_bound_4}. This term can thus be made sufficiently small by reducing $M(\psi)$ and $h$. From estimate \eqref{contractivity}, we then conclude that $\mathcal{J}$ is contractive for sufficiently small $M(\psi)$ and $h$. \\
 \indent On account of Banach's contraction principle, the mapping $\mathcal{J}$ has a unique fixed-point $\psi_h=\mathcal{J}(\psi_h) \in \mathcal{B}_h$, which is, in turn, the unique solution of the nonlinear semi-discrete problem \eqref{Westervelt_ibvp_discr} with approximate initial data \eqref{approx_initial_conditions}.
\end{proof}

\newcommand{\mat}[1]{\mathbf{#1}}
\newcommand{\ten}[1]{\boldsymbol{\mathcal{#1}}}

\section{Computational DG approach for nonlinear sound waves} \label{Sec:Numerical_treatment}

Starting from the semi-discrete equation \eqref{Westervelt_ibvp_discr}, this section describes the numerical treatment and solution process involving the assembly of an equation in matrix-vector form and the time-integration scheme that is used. Computational discontinuous Galerkin approaches for nonlinear acoustic waves that are based on developing a first-order conservative system of equations are investigated in, e.g.,~\cite{kelly2018linear, tripathi2015discontinuous, tripathi2018element}.
\subsection{The matrix form of the semi-discrete problem}
For the purpose of carrying out our numerical experiments, we consider here a more general case than before of having either a non-zero source term or inhomogeneous Dirichlet boundary conditions. We present the numerical treatment of the latter case; the simpler case of having a non-zero source term $f$ can be treated in an analogous manner.\\
\indent Let $\psi = g$ on $\Gamma_{\textup{D}}=\partial\Omega$, where the function $g$ is assumed to be a sufficiently smooth function on $\Gamma_{\textup{D}}$, compatible with initial data.
 The Dirichlet conditions are imposed in a weak sense; see~\cite[Chapter 4, Section 4.2.2]{di2011mathematical} for a detailed explanation. Therefore, in our semi-discrete weak form, the following terms arise additionally on the right-hand side:
\begin{equation*}
 \int_{\Gamma_{\textup{D}}} - c^2\,\tilde{g} \,\left(\nabla v\cdot \vecc{n}\right)~\textup{d}S+\int_{\Gamma_{\textup{D}}} \chi\,\tilde{g} \,v ~\textup{d}S,
\end{equation*}
where, analogously to before, we have used the auxiliary notation $\tilde{g}=g+\tfrac{b}{c^2}\,\dot{g}$. \\
\indent The semi-discrete form of \eqref{Westervelt_ibvp_discr} then reads as
\begin{equation} \label{eq:matrix-equation-reduced}
\mat{M}\vecc{\ddot{\psi}}+ \left(\mat{K}-\mat{D}-\mat{D}^{\top}+\mat{P}\right)(\vecc{\psi} + \tfrac{b}{c^2}\vecc{\dot{\psi}}) - \,\ten{T}[\cdot,\vecc{\dot{\psi}},\vecc{\ddot{\psi}}] = \vecc{w},
\end{equation}
 where $\mat{M}$ denotes the standard mass matrix and $\mat{K}$ the stiffness matrix. In addition, we assemble the nonlinearity tensor $\ten{T}$, the DG penalty matrix $\mat{P}$, the DG jump matrix $\mat{D}=$, and the Dirichlet data vector $\vecc{w}$.\\ 
\indent We therefore have a second-order system of ordinary differential equations with a nonlinear term on the right-hand side, which now remains to be solved by a suitable time-integration scheme. Herein the initial data approximations $(\psi_h(0), \dot{\psi}_{h}(0))=(\psi_{0, h}, \psi_{1, h}) \in V_h \times V_h$ are represented in the finite element basis via the coefficient vectors $\vecc{\psi}_0$ and $\vecc{\psi}_1$ such that for the ODE-system we have $\vecc{\psi}(0)=\vecc{\psi}_0$ and $\vecc{\dot{\psi}}(0)=\vecc{\psi}_1$. 
\subsection{Time integration}
In order to integrate the system of ordinary differential equations \eqref{eq:matrix-equation-reduced} in time, we employ either the Newmark scheme or the Newmark-type Generalized-$\alpha$ method; we refer to~\cite{MKaltenbacher, hoffelner2001finite} for a similar strategy. The nonlinear term is resolved via a fixed-point iteration during the solving stage of the predictor-corrector scheme. The termination criterion that checks the relative change of the solution-vector between iteration steps is employed. In the experiments with realistic physical data, where we observe the nonlinear steepening of the wave front in our computational domain, we choose the Generalized-$\alpha$ scheme because it allows to add targeted numerical damping to the higher modes and subdue Gibbs oscillations.\\
\indent The numerical parameters that are used in the forthcoming experiments can be found in~\Cref{table_param}. The physical and discretization parameters for each experiment are given below in their respective sections. \\
\begin{center}
\renewcommand{\arraystretch}{2.5}
\begin{tabular}{|c||c|c|c|}
	\hline
	                    & Test case 1 & Test case 2 & Test case 3\\
	 \hline	 \hline
	 Newmark scheme & \parbox{2cm}{$\beta_{\textup{nm}}= 0.25$\\ $\gamma_{\textup{nm}}= 0.5$} & \parbox{2cm}{$\beta_{\textup{nm}}=4/9 $\\ $\gamma_{\textup{nm}}= 5/6$} & \parbox{2cm}{$\beta_{\textup{nm}}=4/9 $\\ $\gamma_{\textup{nm}}= 5/6$} \\
	 \hline
	 Generalized-$\alpha$ scheme &  \parbox{2cm}{$\alpha_m= 0$\\ $\alpha_f=0$} & \parbox{2cm}{$\alpha_m= 0$\\ $\alpha_f=1/3$} & \parbox{2cm}{$\alpha_m= 0$\\ $\alpha_f=1/3$} \\
	 \hline
	 Nonlinear iteration & \parbox{2.5cm}{$\textup{TOL}=10^{-5} $\\ $\kappa_{\textup{max}}=100$} & \parbox{2.5cm}{$\textup{TOL}= 10^{-5}$\\ $\kappa_{\textup{max}}=100$} & \parbox{2.5cm}{$\textup{TOL}= 10^{-5}$\\ $\kappa_{\textup{max}}=100$} \\
	 \hline
	 DG penalty & \parbox{2cm}{$\beta=10$} & \parbox{2cm}{$\beta=10$} & \parbox{2cm}{$\beta=250$} \\
	 \hline
\end{tabular}
\captionof{table}{Numerical parameters used in the experiments} \label{table_param} 
\end{center}
~\\[4mm]
In Table~\ref{table_param}, $\beta_{\textup{nm}}$ and $\gamma_{\textup{nm}}$ denote the parameters in the Newmark scheme. The numbers $\alpha_m$ and $\alpha_f$ are the additional parameters that come from the Generalized-$\alpha$ scheme, while $\textup{TOL}$ is the relative tolerance in the termination criterion of the fixed-point-iteration. The number $\kappa_{\textup{max}}$ stands for the maximum number of iterations after which the algorithm should abort. Finally, $\beta$ is the DG-penalty term introduced in \eqref{stabilization}.

				  						  						
\section{Numerical results} \label{Sec:Numerical_results}						  

In this section, we perform numerical experiments to illustrate our theoretical findings. The first two numerical tests are conducted in a two-dimensional computational setting based on a MATLAB implementation. The final, three-dimensional experiment was implemented in SPEED---a parallel, high-order spectral finite-element FORTRAN code~\cite{mazzieri2013speed}.


\subsection{Test case 1: Exact solution known}
\begin{figure}[H]
	\begin{center}
		\hspace*{1cm}\parbox{\textwidth}{
		\input{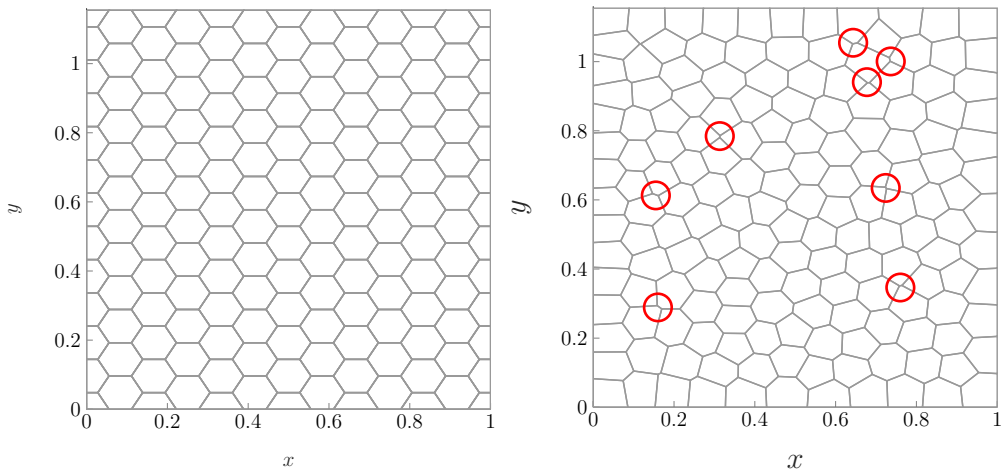}}
		\caption{Computational domain with exemplary polygonal grids. \textbf{(left)} regular hexagonal grid, \textbf{(right)} grid with arbitrary polygonal structure and highlighting of a few exemplary ``small'' edges, making the grid only suitable under the weaker assumptions discussed in~\Cref{rem:remark_trace_inverse_poly}. \label{fig:Ex1_Domain}}
	\end{center}
\end{figure}

\begin{figure}[H]
	\begin{center}
\begin{minipage}{0.49\textwidth}
%
%
\definecolor{mycolor1}{rgb}{1.00000,0.60000,0.00000}%
\begin{tikzpicture}[scale=0.55, font=\Large]

\begin{axis}[%
width=4.568in,
height=3.464in,
at={(0.753in,0.623in)},
scale only axis,
xmode=log,
xmin=0.00872686884847341,
xmax=0.0835169980780651,
xminorticks=true,
xlabel style={at={(0.5, -0.04)}, font=\Large\color{white!15!black}},
xlabel={$1/\sqrt{\textup{No. of elem}}$},
ymode=log,
ymin=8.70505601720726e-08,
ymax=0.0186307995870471,
yminorticks=true,
ylabel style={at={(-0.03, 0.5)}, font=\Large\color{white!15!black}},
ylabel={$L^2$ error},
axis background/.style={fill=white},
xmajorgrids,
ymajorgrids,
legend style={at={(0.02,0.98)}, anchor=north west, legend cell align=left, align=left, draw=white!15!black}
]
\addplot [color=blue, line width=2.0pt]
table[row sep=crcr]{%
	0.0145447814141224	9.23087460172776e-06\\
	0.0596549986271894	0.000636886106340099\\
};
\addlegendentry{$h^{3}$}

\addplot [color=blue, dashdotted, line width=2.0pt]
table[row sep=crcr]{%
	0.0145447814141224	1.34261053343304e-07\\
	0.0596549986271894	3.79934397993947e-05\\
};
\addlegendentry{$h^{4}$}

\addplot [color=mycolor1, mark=diamond*, mark options={solid, fill=mycolor1, mycolor1}, line width=2.0pt]
table[row sep=crcr]{%
	0.0596549986271894	0.00242974447157846\\
	0.0347734071210361	0.000341024528757495\\
	0.0233890134862494	8.8281119026889e-05\\
	0.0182635074342946	3.92934046411971e-05\\
	0.0145447814141224	1.8981662209964e-05\\
};
\addlegendentry{$p=2:$ (reg.)}

\addplot [color=mycolor1, mark=square*, mark options={solid, fill=mycolor1, mycolor1}, line width=2.0pt]
table[row sep=crcr]{%
	0.0596549986271894	0.000140077638624087\\
	0.0347734071210361	1.46243448633194e-05\\
	0.0233890134862494	2.97187496602835e-06\\
	0.0182635074342946	1.2492882179431e-06\\
	0.0145447814141224	4.35252800860363e-07\\
};
\addlegendentry{$p=3:$ (reg.)}

\addplot [color=black!20!green, mark=diamond*, mark options={solid, fill=black!20!green, black!20!green}, line width=2.0pt]
table[row sep=crcr]{%
	0.0596549986271894	0.00273982346868339\\
	0.0347734071210361	0.000473324075292617\\
	0.0233890134862494	0.000118290502007015\\
	0.0182635074342946	5.05456836313172e-05\\
	0.0145447814141224	2.38243856984137e-05\\
};
\addlegendentry{$p=2:$ (rand.)}

\addplot [color=black!20!green, mark=square*, mark options={solid, fill=black!20!green, black!20!green}, line width=2.0pt]
table[row sep=crcr]{%
	0.0596549986271894	0.000110180054995583\\
	0.0347734071210361	1.20490810843636e-05\\
	0.0233890134862494	2.5642888528406e-06\\
	0.0182635074342946	1.12347947292158e-06\\
	0.0145447814141224	3.83052681826157e-07\\
};
\addlegendentry{$p=3:$ (rand.)}

\end{axis}
\end{tikzpicture}%
\end{minipage}%
\begin{minipage}{0.49\textwidth}
\hspace*{1mm}			
%
%
\definecolor{mycolor1}{rgb}{1.00000,0.60000,0.00000}%
\begin{tikzpicture}[scale=0.55, font=\Large]

\begin{axis}[%
width=4.568in,
height=3.464in,
at={(0.766in,0.626in)},
scale only axis,
xmode=log,
xmin=0.00872686884847341,
xmax=0.0835169980780651,
xminorticks=true,
xlabel style={at={(0.5, -0.04)}, font=\Large\color{white!15!black}},
xlabel={$1/\sqrt{\textup{No. of elem}}$},
ymode=log,
ymin=4.01772518063185e-05,
ymax=1.94041032719861,
yminorticks=true,
ylabel style={at={(-0.03, 0.5)}, font=\Large\color{white!15!black}},
ylabel={$|\nabla_h (\cdot)|_{L^2}$ error},
axis background/.style={fill=white},
xmajorgrids,
ymajorgrids,
legend style={at={(0.02,0.98)}, anchor=north west, legend cell align=left, align=left, draw=white!15!black}
]
\addplot [color=blue, line width=2.0pt]
table[row sep=crcr]{%
	0.0145447814141224	0.00380791199492278\\
	0.0596549986271894	0.0640569395017794\\
};
\addlegendentry{$h^2$}

\addplot [color=blue, dashdotted, line width=2.0pt]
table[row sep=crcr]{%
	0.0145447814141224	5.53852476103665e-05\\
	0.0596549986271894	0.0038213166380406\\
};
\addlegendentry{$h^3$}

\addplot [color=mycolor1, mark=diamond*, mark options={solid, fill=mycolor1, mycolor1}, line width=2.0pt]
table[row sep=crcr]{%
	0.0596549986271894	0.159916844290133\\
	0.034773407121036	0.0519372855566819\\
	0.0233890134862494	0.0231460379813486\\
	0.0182635074342946	0.0140339159715621\\
	0.0145447814141224	0.00886902299030062\\
};
\addlegendentry{$p=2$ (reg.)}

\addplot [color=mycolor1, mark=square*, mark options={solid, fill=mycolor1, mycolor1}, line width=2.0pt]
table[row sep=crcr]{%
	0.0596549986271894	0.0154005117980149\\
	0.034773407121036	0.00287337139557262\\
	0.0233890134862494	0.00085210173517689\\
	0.0182635074342946	0.000401186517056585\\
	0.0145447814141224	0.000200886259031593\\
};
\addlegendentry{$p=3$ (reg.)}

\addplot [color=black!20!green, mark=diamond*, mark options={solid, fill=black!20!green, black!20!green}, line width=2.0pt]
table[row sep=crcr]{%
	0.0596549986271894	0.161700860599884\\
	0.034773407121036	0.0556490197654581\\
	0.0233890134862494	0.0251785457501266\\
	0.0182635074342946	0.0153468810742043\\
	0.0145447814141224	0.00969505047909195\\
};
\addlegendentry{$p=2$ (rand.)}

\addplot [color=black!20!green, mark=square*, mark options={solid, fill=black!20!green, black!20!green}, line width=2.0pt]
table[row sep=crcr]{%
	0.0596549986271894	0.013769213413617\\
	0.034773407121036	0.00266175412587811\\
	0.0233890134862494	0.000836760858214949\\
	0.0182635074342946	0.000399547766666226\\
	0.0145447814141224	0.000199993512649204\\
};
\addlegendentry{$p=3$ (rand.)}

\end{axis}
\end{tikzpicture}%
\end{minipage}
	\end{center}
		\caption{$L^2$ and $|\nabla_h (\cdot)|_{L^2}$ errors of $\psi_h$ at final time for four sequentially refined polygonal meshes and second- and third-order polynomials, comparing a sequence of regular polygonal grids (reg.) with a sequence of irregular ones (rand.).
	\label{fig:Ex1_Error_h}}
\end{figure}
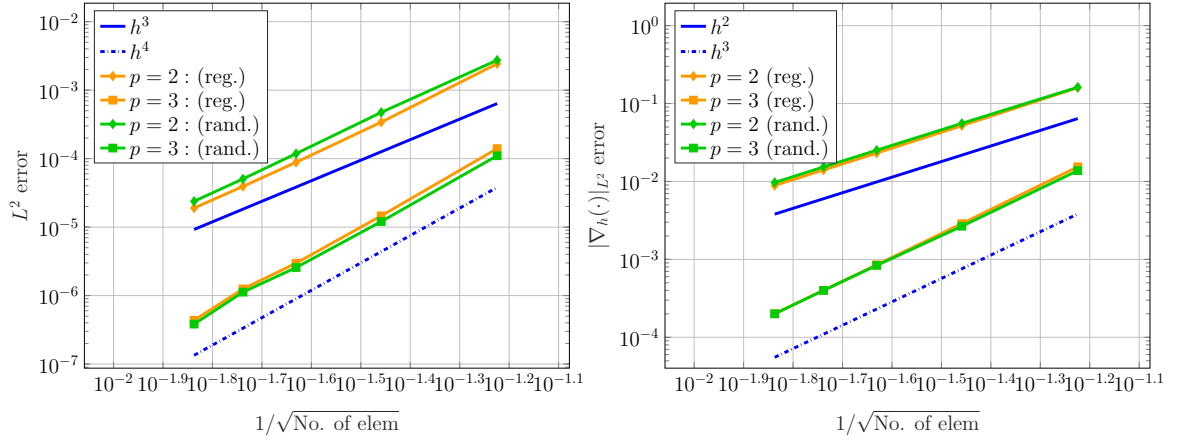

In our first example, we simulate the Westervelt equation \eqref{Westervelt_potential} with a given source term $f$ on the right hand side, which we choose as
\begin{align*}
	f &= \left[16\pi^2(c^2-1)\sin(4\pi t)+64\pi^3b\cos(4\pi t)\right] \sin(4\pi x) \\
	&\hspace*{1cm}+ \left[64\pi^3 k\sin(4\pi t)\cos(4\pi t)\right]\sin(4\pi x)^2.
\end{align*}
In this way, the exact solution is given by $\psi = \sin(4\pi x)\sin(4\pi t)$, which we use in the error analysis. In this numerical experiment, all the physical quantities involved are assumed to be dimensionless. Our computational domain is given by the rectangle $\Omega = (0, 1)\times (0, \frac{2}{3}\sqrt{3})$. We tessellate it with $N_{\textup{elem}}$ polygonal elements in two ways: a regular hexagonal pattern and a random way using polygons with different number of edges each; see Figure~\ref{fig:Ex1_Domain} for exemplary depictions of the resulting grids. The initial conditions and Dirichlet boundary data are set to correspond to the values of $\psi$ at time zero and on the boundary, respectively.\\
\indent  We choose the coefficents in the equation to be $c=1$, $b=10^{-5}$, $\beta_a = 10^{-4}$, and the mass density is $\rho=1$. The time-discretization is conducted with final time $T = 0.8$ and the Newmark scheme, where the time stepsize is always adapted in such a way, that the time-discretization error does not dominate in~\Cref{fig:Ex1_Error_h} and the convergence with respect to the number of elements can be observed.


We perform this experiment on a sequence of regular grids that satisfy our mesh assumptions and on an unstructured grid which might not; see~\Cref{fig:Ex1_Domain}. Such unstructured meshes satisfy the weaker assumptions discussed in~\Cref{rem:remark_trace_inverse_poly}, which is why we want to test if we still observe the same order of convergence when using them.\\
\indent \Cref{fig:Ex1_Error_h} displays convergence results for five sequentially refined polygonal meshes, where we have employed polynomials of degree $p=2$ and $p=3$. As a reference, on the five levels the unstructured grid consists out of 281, 827, 1828, 2998, and 4727 elements, respectively. As expected in practice for DG methods, the $L^2$-error of the acoustic velocity potential converges with the order $h^{p+1}$; see, e.g.,~\cite{grote2006discontinuous}. 




\subsection{Test case 2: Exact solution unknown}\label{subsec:Ex2}

\begin{figure}[H]
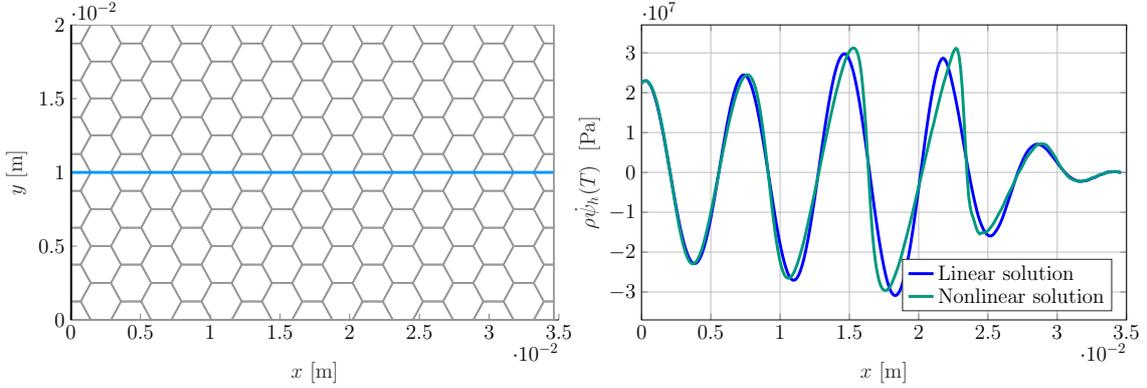

	\begin{center}
		\input{images/Ex2_Domain.tex}\hspace*{-1mm}\input{images/Ex2_Wave.tex}
		\caption{\textbf{(left)} Computational domain with exemplary polygonal grid. The left boundary is the excitation part, where the Dirichlet data $g$ is non-zero. At remaining boundary parts homogeneous Dirichlet data are imposed. The blue line is the axis of symmetry over which the solution is evaluated on the right. \textbf{(right)} Computed acoustic pressure field $u_h=\rho\dot{\psi_h}$ over the horizontal axis of symmetry of the channel at final time $T$ of a simulation with increased amplitude $A=\SI{0.0175}{m^2/s^2}$, plotted together with the pressure wave obtain by solving the linear damped wave equation with same boundary and initial data.
			\label{fig:Ex2_Setup}}
	\end{center}
\end{figure}

Our second example features a more realistic setting. The computational domain is chosen to be a rectangle with dimensions $H=\SI{0.02}{m}$ and $L=\frac{3}{\sqrt{3}}\cdot 0.02\,\textup{m}$; see Figure~\ref{fig:Ex2_Setup} left. The physical parameters are now set to $$c=\SI{1500}{m/s},\ b=\SI{6e-9}{m^2/s},\ \beta_a = 7,\ \text{and} \ \rho=1000\,\textup{kg}/\textup{m}^3.$$ The time horizon is $T=\SI{2.4e-5}{s}$, resolved by a step size of $\textup{dt}=\SI{2e-9}{s}$. Instead of a non-zero source term $f$, we employ inhomogeneous Dirichlet conditions. The excitation signal is given in the form $g(x,y, t) = g^{(s)}(x,y) \cdot g^{(t)}(t)$. Herein the temporal part responsible for the initialization of the wave oscillations is given by

\begin{equation*}\label{eq:Excitation}
	g^{(t)}(t)=\begin{cases}
		\left(\frac{ft}{2}\right)^2\,A\,\sin(\omega t), & t<\frac{2}{f}\,\textup{s}, \\
		\hspace*{1.15cm}A\,\sin(\omega t), & t \geq \frac{2}{f}\,\textup{s},
	\end{cases}
\end{equation*}
while the spatial part is given by a mollifier-type function in order to get a spatially smooth transition between the inhomogeneous excitation and the homogeneous remaining boundary data. In particular, we have

\begin{equation*}
	g^{(s)}(x,y)=\begin{cases}
		0, &x=0\,\textup{m}, y=0\,\textup{m}\\
		\exp\left(1-\frac{1}{1-\left|\frac{1}{0.005}(y-0.005)\right|^2}\right), & x=0\,\textup{m}, 0\,\textup{m}<y<0.005\,\textup{m}\\
		1, &x=0\,\textup{m}, 0.005\,\textup{m}\leq y\leq 0.015\,\textup{m}\\
		\exp\left(1-\frac{1}{1-\left|\frac{1}{0.005}(y-0.015)\right|^2}\right), & x=0\,\textup{m}, 0.015\,\textup{m}<y<0.02\,\textup{m}\\
		0, & x=0\,\textup{m}, y=0.02\,\textup{m}\\
		0, & 0\,\textup{m}<x
	\end{cases}
\end{equation*}

We therefore employ a spatially smooth, temporally modulated sinusoidal excitation with driving frequency $f=\SI{210}{\kilo \hertz}$ and amplitude $A=\SI{0.01}{m^2/s^2}$, where $\omega$ denotes the angular frequency $\omega = 2\pi f$. \Cref{fig:Ex2_Setup} displays a plot of the acoustic pressure $u_h=\rho\dot{\psi}_h$ along the horizontal axis of symmetry of the channel at final time $T$. To better observe the nonlinear steepening, the figure also contains a pressure wave obtained by solving the linear damped wave equation (i.e., Westervelt's equation with $k=0$).


We next want to analyze the behavior of the numerical solution with respect to $h$- and $p$-refinement. However, in this more realistic setting an exact solution is unknown. Therefore, instead of tracking the deviation from a given solution, we track a given quantity of interest. Here we choose to compute 
\begin{align} \label{QoI}
Q(\psi_h):=\|{\psi}_h\|_{L^{\infty}(0,T; L^2(\Omega))},
\end{align}
on different discretization levels. We note that  $$Q(\psi)- \|{\psi}_h-{\psi}\|_{L^{\infty}(0,T; L^2(\Omega))} \leq Q(\psi_h)\leq Q(\psi) + \|{\psi}_h-{\psi}\|_{L^{\infty}(0,T; L^2(\Omega))},$$
and so we expect that, for $p$ fixed, $Q(\psi_h)$ behaves asymptotically as $q_1+q_2 \cdot h^{p+1}$ for some constants $q_1$ and $q_2$.\\

\begin{figure}[H]
	\begin{center}
%
%
\definecolor{mycolor1}{rgb}{1.00000,0.60000,0.00000}%
\begin{tikzpicture}[scale = 0.55, font=\Large]

\begin{axis}[%
width=4.514in,
height=3.464in,
at={(0.916in,0.576in)},
scale only axis,
xmin=0,
xmax=0.00325,
xlabel style={at={(0.5,-0.025)}, font=\Large\color{white!15!black}},
xlabel={$h$},
ymin=0.000132225,
ymax=0.00013228,
ytick={0.000132225,0.000132235,0.000132245,0.000132255,0.000132265,0.000132275},
y tick label style={/pgf/number format/.cd, precision=5, /tikz/.cd},
ylabel style={font=\Large\color{white!15!black}},
axis background/.style={fill=white},
xmajorgrids,
ymajorgrids,
legend style={at={(0.03,0.03)}, anchor=south west, legend cell align=left, align=left, draw=white!15!black}
]
\addplot [color=blue, line width=2.0pt]
table[row sep=crcr]{%
	0	0.000132273477136618\\
	0.000432027338244996	0.000132273459557213\\
	0.000586322816189637	0.000132273417500905\\
	0.00070975919854535	0.000132273349079609\\
	0.000802336485312135	0.000132273268021514\\
	0.000894913772078919	0.000132273153480366\\
	0.000987491058845704	0.000132272997301618\\
	0.00104920925002356	0.000132272865620825\\
	0.00111092744120142	0.000132272708533729\\
	0.00117264563237927	0.000132272522965216\\
	0.00123436382355713	0.00013227230566445\\
	0.00129608201473499	0.000132272053204877\\
	0.00135780020591284	0.000132271761984218\\
	0.0014195183970907	0.000132271428224476\\
	0.00148123658826856	0.000132271047971932\\
	0.00154295477944641	0.000132270617097147\\
	0.00160467297062427	0.000132270131294961\\
	0.00166639116180213	0.000132269586084493\\
	0.00172810935297998	0.00013226897680914\\
	0.00175896844856891	0.000132268646643896\\
	0.00178982754415784	0.000132268298636581\\
	0.00182068663974677	0.000132267932155699\\
	0.0018515457353357	0.000132267546558771\\
	0.00188240483092462	0.000132267141192335\\
	0.00191326392651355	0.000132266715391947\\
	0.00194412302210248	0.000132266268482179\\
	0.00197498211769141	0.000132265799776623\\
	0.00200584121328034	0.000132265308577886\\
	0.00203670030886926	0.000132264794177593\\
	0.00206755940445819	0.000132264255856388\\
	0.00209841850004712	0.00013226369288393\\
	0.00212927759563605	0.000132263104518898\\
	0.00216013669122498	0.000132262490008988\\
	0.00219099578681391	0.00013226184859091\\
	0.00222185488240283	0.000132261179490395\\
	0.00225271397799176	0.000132260481922193\\
	0.00228357307358069	0.000132259755090066\\
	0.00231443216916962	0.000132258998186798\\
	0.00234529126475855	0.000132258210394188\\
	0.00237615036034748	0.000132257390883054\\
	0.0024070094559364	0.000132256538813231\\
	0.00243786855152533	0.000132255653333571\\
	0.00246872764711426	0.000132254733581943\\
	0.00249958674270319	0.000132253778685235\\
	0.00253044583829212	0.000132252787759351\\
	0.00256130493388104	0.000132251759909214\\
	0.00259216402946997	0.000132250694228763\\
	0.0026230231250589	0.000132249589800955\\
	0.00265388222064783	0.000132248445697764\\
	0.00268474131623676	0.000132247260980182\\
	0.00271560041182569	0.000132246034698219\\
	0.00274645950741461	0.0001322447658909\\
	0.00277731860300354	0.000132243453586271\\
	0.00280817769859247	0.000132242096801393\\
	0.0028390367941814	0.000132240694542345\\
	0.00286989588977033	0.000132239245804223\\
	0.00290075498535926	0.000132237749571142\\
	0.00293161408094818	0.000132236204816233\\
	0.00296247317653711	0.000132234610501644\\
	0.00299333227212604	0.000132232965578543\\
	0.00302419136771497	0.000132231268987113\\
	0.0030550504633039	0.000132229519656555\\
};
\addlegendentry{$1.322735\cdot 10^{-4}-504.613929\cdot h^{4}$}

\addplot [color=mycolor1, mark=diamond*, mark options={solid, fill=mycolor1, mycolor1},  line width=2.0pt]
table[row sep=crcr]{%
	0.0030550504633039	0.000132230258843918\\
	0.00254587538608658	0.000132250470425968\\
	0.00218217890235994	0.000132261825860311\\
	0.00190940653956494	0.00013226765779297\\
	0.00169725025739106	0.000132269729329648\\
	0.00152752523165196	0.000132270841239314\\
	0.00132828281013214	0.000132272048084694\\
	0.00105346567700135	0.000132272795715462\\
	0.000848625128695533	0.000132272972215603\\
};
\addlegendentry{$\|{\psi_h}\|_{L^{\infty}(0,T;L^2(\Omega))}$}

\end{axis}
\end{tikzpicture}%
%
%
\definecolor{mycolor1}{rgb}{1.00000,0.60000,0.00000}%
\pgfplotsset{
	log x ticks with fixed point/.style={
		xticklabel={
			\pgfkeys{/pgf/fpu=true}
			\pgfmathparse{exp(\tick)}%
			\pgfmathprintnumber[fixed relative, precision=3]{\pgfmathresult}
			\pgfkeys{/pgf/fpu=false}
		}
	},
	log y ticks with fixed point/.style={
		yticklabel={
			\pgfkeys{/pgf/fpu=true}
			\pgfmathparse{exp(\tick)}%
			\pgfmathprintnumber[fixed relative, precision=3]{\pgfmathresult}
			\pgfkeys{/pgf/fpu=false}
		}
	}
}
\begin{tikzpicture}[scale = 0.55, font=\Large]

\begin{axis}[%
width=4.514in,
height=3.464in,
at={(0.916in,0.576in)},
scale only axis,
log x ticks with fixed point,
xmode=log,
xmin=1,
xmax=7,
xminorticks=true,
xlabel style={at={(0.5,-0.025)},font=\Large\color{white!15!black}},
xlabel={$p$},
xtick={1,2,3,4,5,6,7},
ymode=log,
ymin=1e-11,
ymax=1e-05,
yminorticks=true,
ylabel style={font=\Large\color{white!15!black}},
ylabel={},
axis background/.style={fill=white},
xmajorgrids,
xminorgrids,
ymajorgrids,
yminorgrids,
legend style={at={(0.03,0.03)}, anchor=south west, legend cell align=left, align=left, draw=white!15!black}
]
\addplot [color=mycolor1, mark=diamond*, mark options={solid, fill=mycolor1, mycolor1}, line width=2.0pt]
  table[row sep=crcr]{%
1	1.36894895496668e-06\\
2	1.49893289888038e-07\\
3	2.26295740321156e-08\\
4	1.85270870211229e-09\\
5	2.21753781786891e-10\\
6	1.50650514280223e-10\\
7	4.02252866974796e-11\\
};
\addlegendentry{$\|{\psi}_h\|_{L^{\infty}(0,T;L^2(\Omega))}-\textup{reference value}$}

\end{axis}
\end{tikzpicture}%
		\caption{\textbf{(left)} The quantity of interest $Q(\psi_h)=\|{\psi}_h\|_{L^{\infty}(0,T;L^2(\Omega))}$ computed on a sequence of polygonal meshes with third-order polynomials and a least-square fitted extrapolation-curve of order $h^4$. \textbf{(right)} Deviation of the computed quantity of interest $Q(\psi_h)$ from a reference value as a function of $p$ on a fixed mesh.
		\label{fig:Ex2_res}}
	\end{center}
\end{figure}
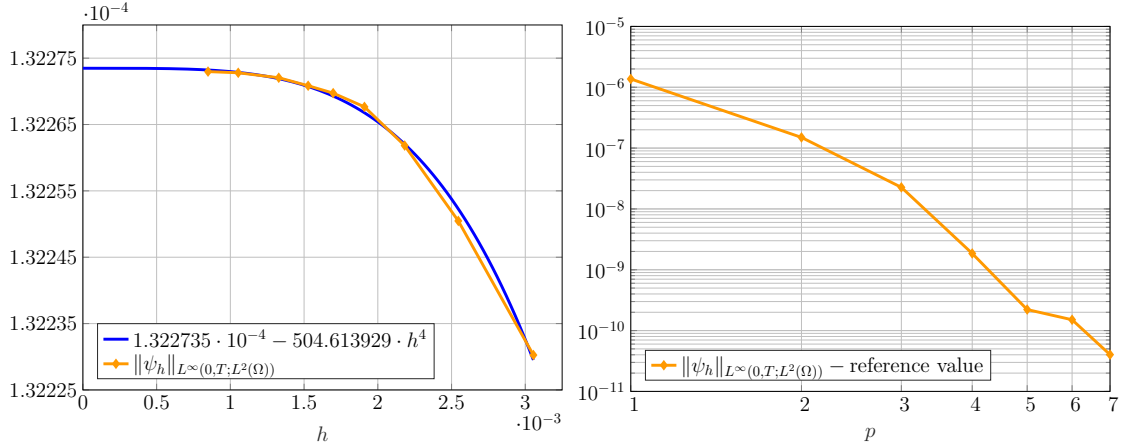

\paragraph{\textbf{\mathversion{bold} $h$-refinement}} We restrict ourselves to structured, quasi-uniform polygonal meshes consisting of $N_{\textup{elem}}\in\lbrace 220, 312, 420, 544, 684, 840, 1104, 1740,$ $2664\rbrace$ elements. The polynomial degree is set to $p=3$. Values $Q(\psi_h)$ for these levels of refinement are plotted in~\Cref{fig:Ex2_res} on the left. In order to observe the convergence order, we perform a least-square fit of the $(h,Q(\psi_h))$ data pairs. We obtain a fitted curve  $$Q_f(h)=q_1+q_2\cdot h^{4},$$ where $q_1$ and $q_2$ are subject to the least-square fit. The fitted curve $Q_f$ with optimized parameters reads approximately as $$Q_f(h)=1.322735\cdot 10^{-4}-504.613929\cdot h^{4}$$ and is plotted in~\Cref{fig:Ex2_res} on the left as well. As expected, we observe $O(h^{p+1})$ convergence. The extrapolated value for the quantity of interest $Q_f(0)$ evaluates to $1.322735\cdot 10^{-4}$.\\


\paragraph{\textbf{\mathversion{bold} $p$-refinement}} We also perform a refinement analysis of the quantity of interest with respect to the polynomial degree. Here we choose a fixed mesh with $N_{\textup{elem}}=312$ elements and successively increase the polynomial degree $p\in\lbrace 1,2,...,7\rbrace$, where it should be mentioned that even though our theory holds for $p\geq 2$, the case $p=1$ yields a similar result as well. The deviations of the resulting quantities from a reference value are plotted in~\Cref{fig:Ex2_res} on the right. As the reference, we choose the extrapolated value of the $h$-refinement with degree $p=4$. 
For $p=7$, the quantity of interest evaluates to $1.322730\cdot 10^{-4}$ with the deviation from the extrapolated reference value below $7.6\cdot 10^{-5}\,\%$. 

\subsection{Test case 3: A three-dimensional example}

\begin{figure}[H]
	\begin{center}
		\includegraphics[trim=0mm 0mm 0mm 0mm, clip, scale=0.5]{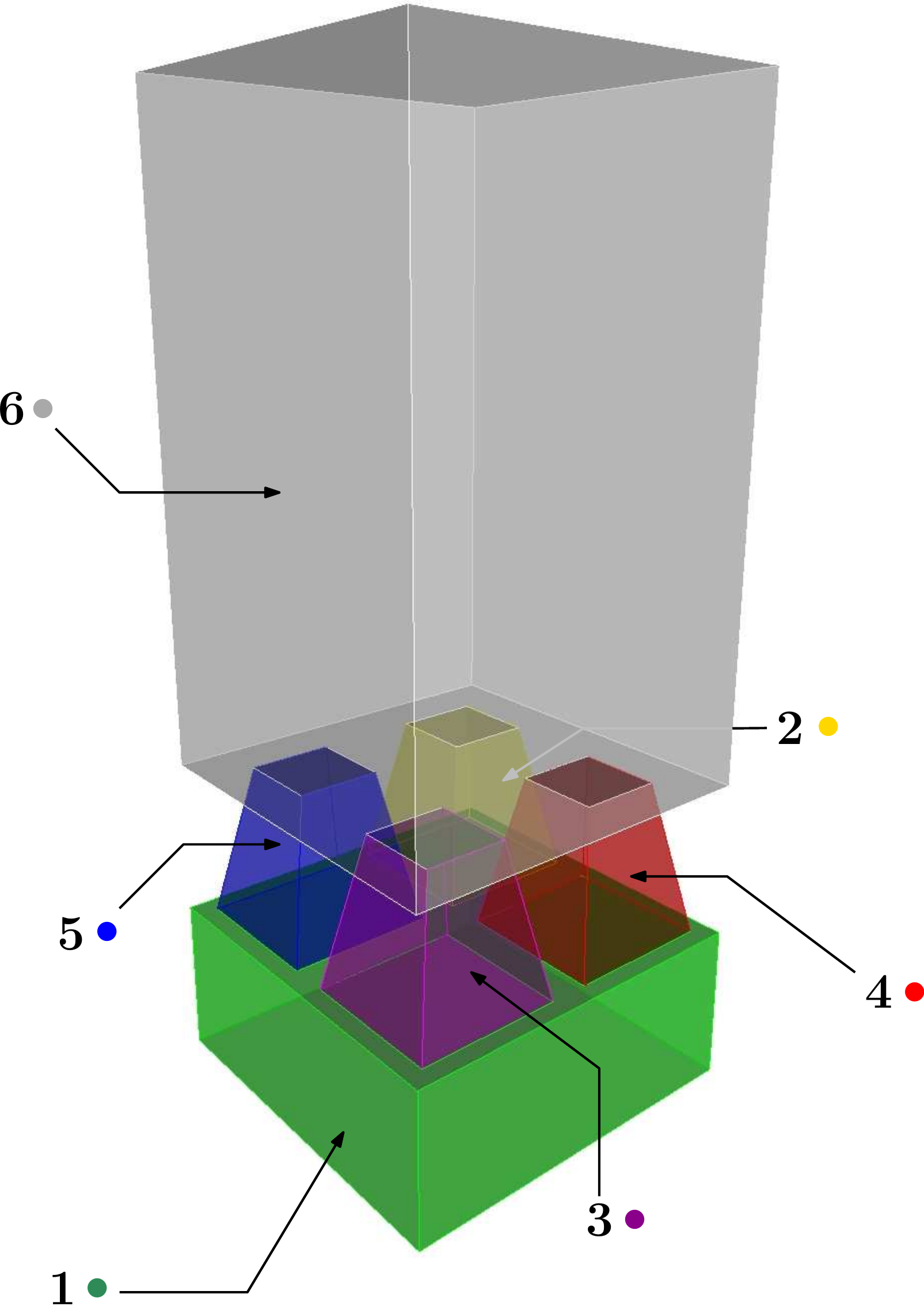}
		\caption{Test case 3: Computational domain with six different material blocks. Each block has its own material parameters, mesh, and ansatz space. Interface coupling is done via the DG approach. 
			\label{fig:Ex3_Domain}}
	\end{center}
\end{figure}

Our final example is performed a three dimensional setting, where we use the discontinuous Galerkin approach in a hybrid way. Figure \ref{fig:Ex3_Domain} shows an image of the computational geometry consisting of six different material blocks. The blocks differ in all relevant material parameters, given in Table \ref{tab:phys_param}. We therefore in this experiment solve Westervelt's equation
	\begin{equation*} \label{Westervelt_potential_jumps}
\begin{aligned}
(1-2k \dot{\psi})\ddot{\psi}-\textup{div}(c^2 \nabla \psi) - \textup{div}(b \nabla \dot{\psi})=0
\end{aligned}
\end{equation*}
with coefficients in $L^\infty(\Omega)$. Now the stabilization function has the form
\begin{equation*} \label{stabilization_blocks}
\begin{aligned}
\chi_{\vert F}= \begin{cases}
 \, \beta \, \displaystyle \max_{\kappa \in \{\kappa^+, \kappa^-\}} c^2_{\vert \kappa} \frac{p^2}{h_\kappa} \quad \quad &\text{for all  } F \in \mathcal{F}_h^i, \quad F \subset \partial \kappa^+ \cap \partial \kappa^-, \\[5mm]
\ \beta c^2_{\vert \kappa} \, \dfrac{p^2}{h_\kappa}  &\text{for all  } F \in \mathcal{F}_h^b, \quad F \subset \partial \kappa.
\end{cases}
\end{aligned}
\end{equation*}
To save computational power, the blocks are meshed and discretized individually via conforming spectral elements within each block, while the discontinuous Galerkin approach deals with the non-matching grids on the interfaces.

\begin{center}
	\renewcommand{\arraystretch}{1.5}
	\begin{tabular}{|c||c|c|c|c|c|c|}
		\hline
		Material block & 1 $\color{sgreen}\bullet$ & 2 $\color{syellow}\bullet$ & 3 $\color{spink}\bullet$ & 4 $\color{sred}\bullet$ & 5  $\color{sblue}\bullet$& 6 $\color{sgray}\bullet$ \\
		\hline	 \hline
		 $c\,\left[\si{m/s}\right]$ & 1500& 1000& 1000& 3000&750 &1500 \\
		\hline
		 $b\,\left[\si{m^2/s}\right]$ & $\num{6d-9}$& $\num{4d-9}$& $\num{4d-2}$&$\num{4d-9}$ & $\num{6d-9}$& $\num{6d-9}$\\
		\hline
		 $\beta_a$ & 5& 4& 4& 7& 4& 5\\
		\hline
		 $\rho\,\left[\si{kg/m^3}\right]$ & 1000 & 1250& 1250 & 2000 & 1500 & 1000 \\
		\hline
	\end{tabular}
	\captionof{table}{Physical parameters of different material blocks.} \label{tab:phys_param} 
\end{center}~\\

\begin{figure}[H]
	\begin{center}
		\parbox{22cm}{\includegraphics[trim=16cm 0mm 22cm 0mm, clip, scale=0.175]{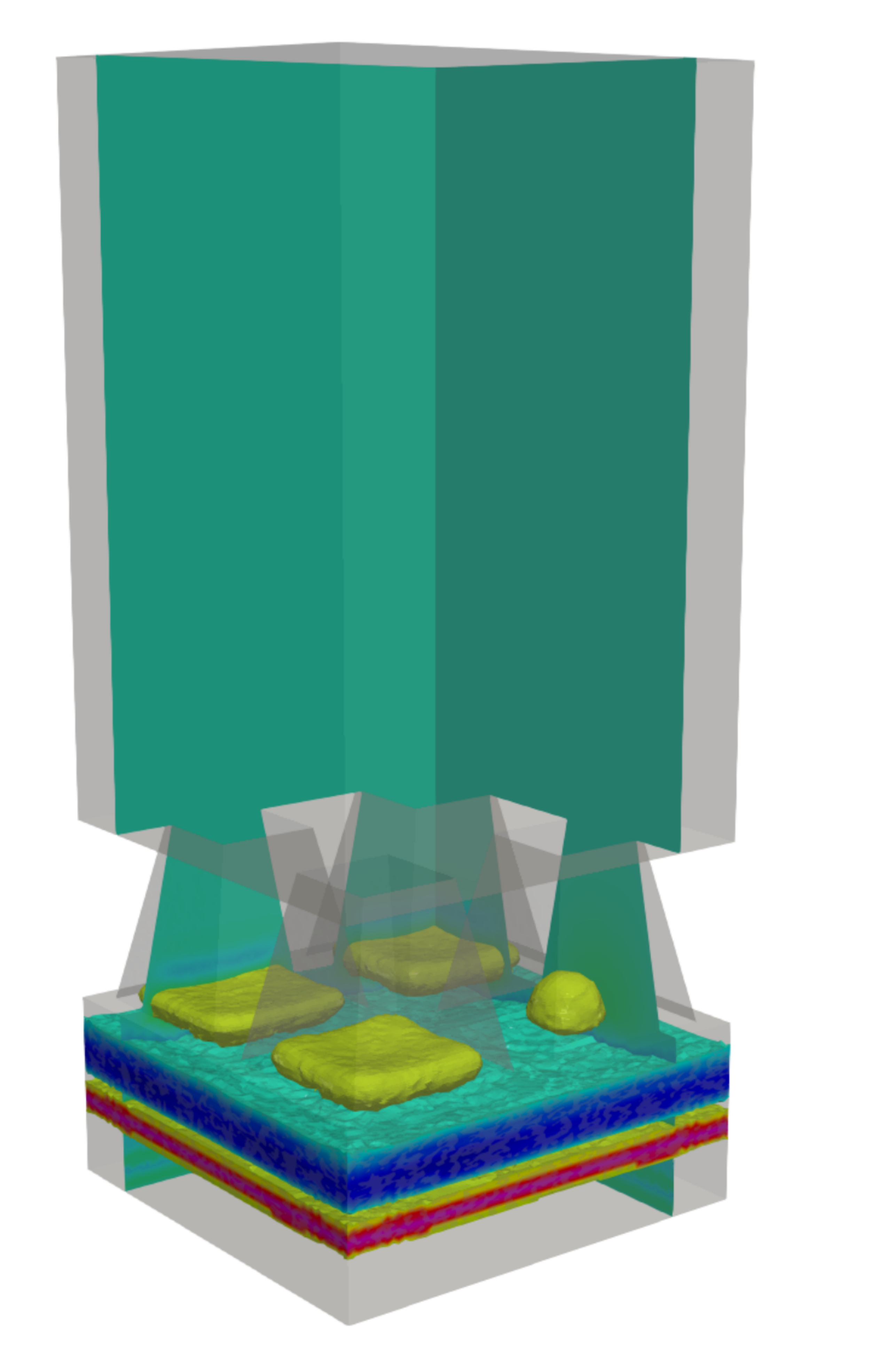}		\includegraphics[trim=16cm 0mm 22cm 0mm, clip, scale=0.175]{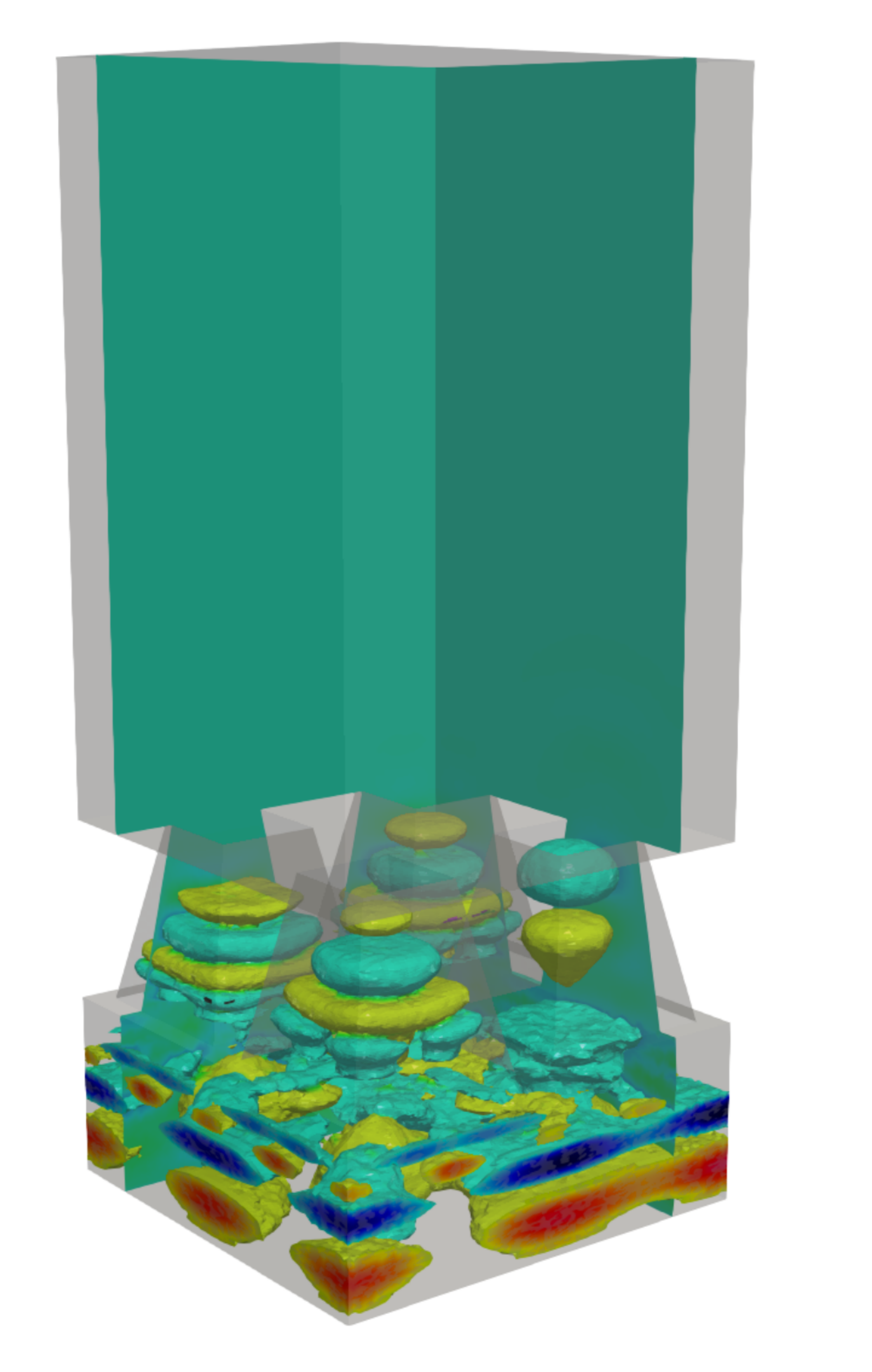}
		\includegraphics[trim=16cm 0mm 22cm 0mm, clip, scale=0.175]{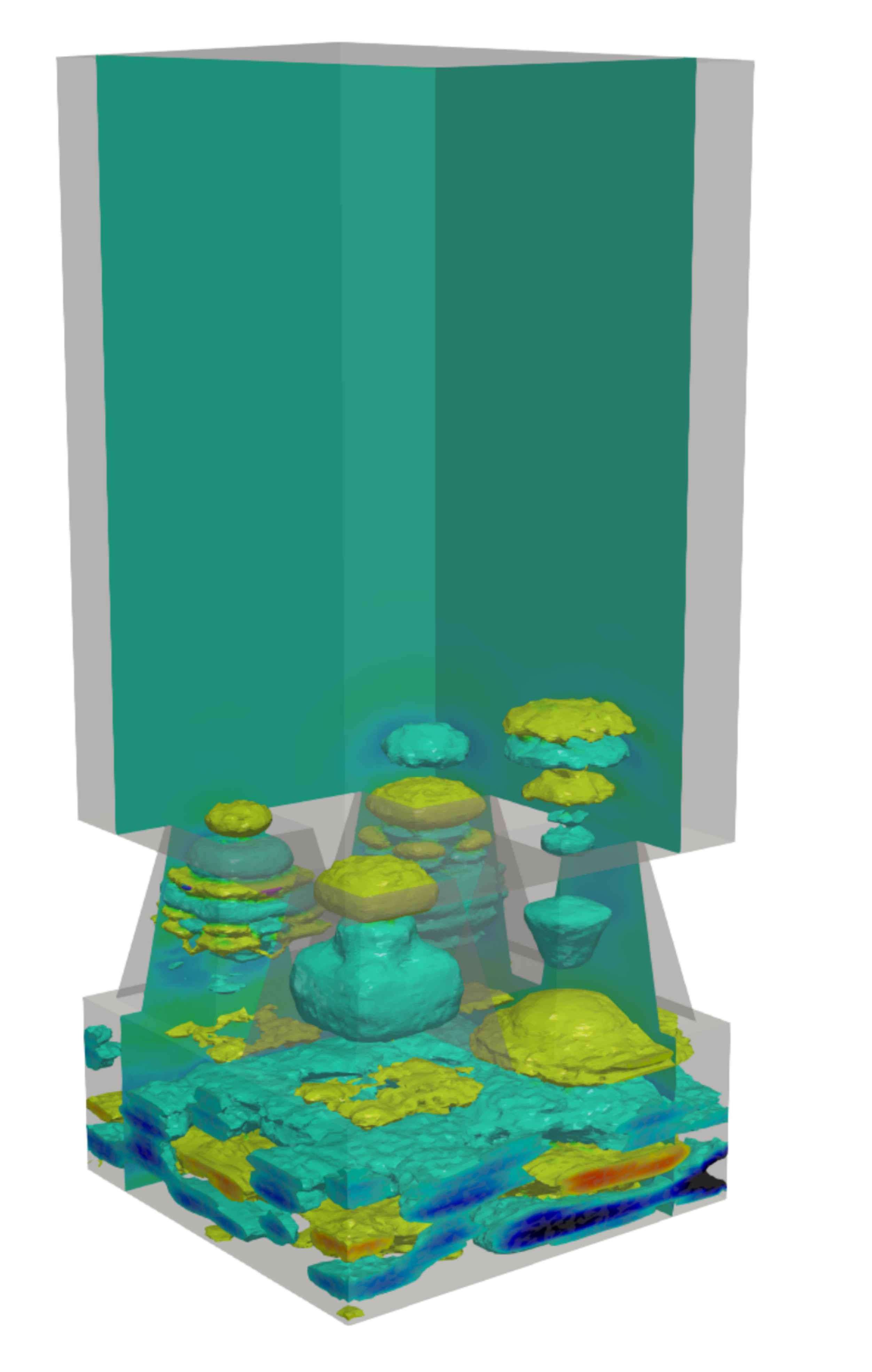}
		\raisebox{1.25cm}{\includegraphics[trim=0cm 0mm 0cm 0mm, clip, scale=0.3]{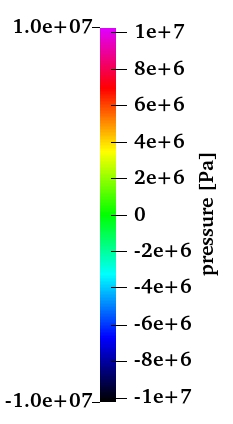}}}
		\caption{Iso-volumina at time steps \textbf{(left)} 12000 \textbf{(middle)} 19000 \textbf{(right)} 24600 of the highest acoustic pressure amplitudes (in absolute value) $u_h=\rho\dot{\psi}_h$ during wave propagation through the four connecting channels of the computational domain. The orientation of the images is the same as in \Cref{fig:Ex3_Domain}.\label{fig:Ex3_pics}}
	\end{center}
\end{figure}

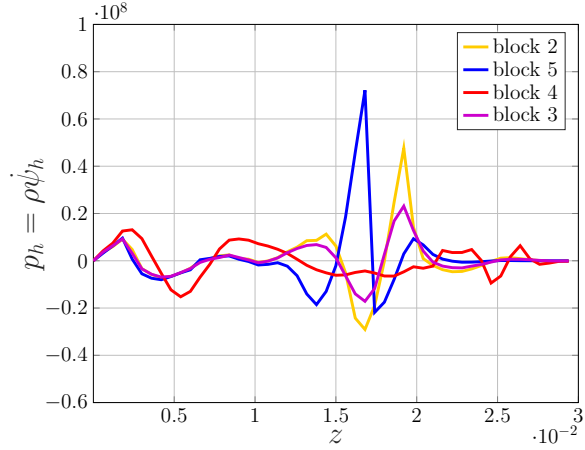
\begin{figure}[H]
	\begin{center}
%
%
\definecolor{mycolor1}{rgb}{1.00000,0.80000,0.00000}%
\definecolor{mycolor2}{rgb}{0.80000,0.00000,0.80000}%
\begin{tikzpicture}[scale = 0.55, font=\Large]

\begin{axis}[%
width=4.568in,
height=3.589in,
at={(0.766in,0.5in)},
scale only axis,
xmin=0,
xmax=0.03,
xlabel style={font=\Large\color{white!15!black}},
xlabel={\huge $z$},
xtick={0.005, 0.01, 0.015, 0.02, 0.025, 0.03},
ymin=-60000000,
ymax=100000000,
ylabel style={at={(-0.03,0.5)}, font=\Large\color{white!15!black}},
ylabel={\huge $p_h=\rho\dot{\psi}_h$},
axis background/.style={fill=white},
title style={font=\bfseries},
xmajorgrids,
ymajorgrids,
legend style={at={(0.75,0.85)}, anchor=west, legend cell align=left, align=left, draw=white!15!black}
]
\addplot [color=mycolor1, line width=2.0pt]
  table[row sep=crcr]{%
0	0\\
0.000600002706050873	3575800\\
0.00119999796152115	5979600\\
0.00180000066757202	9154700\\
0.00240000337362289	4773300\\
0.00299999862909317	-3318400\\
0.00360000133514404	-5802100\\
0.00419999659061432	-7120100\\
0.00479999929666519	-6597600\\
0.00540000200271606	-5139000\\
0.00599999725818634	-3157100\\
0.00659999996423721	-173490\\
0.00720000267028809	848990\\
0.00779999792575836	2070000\\
0.00840000063180923	2489300\\
0.00900000333786011	920960\\
0.00959999859333038	670260\\
0.0102000012993813	-672020\\
0.0107999965548515	-66838\\
0.0113999992609024	1312900\\
0.0120000019669533	3819400\\
0.0125999972224236	5760500\\
0.0131999999284744	8496500\\
0.0138000026345253	8632900\\
0.0143999978899956	11273000\\
0.0150000005960464	6242000\\
0.0156000033020973	-5107600\\
0.0161999985575676	-24269000\\
0.0168000012636185	-29073000\\
0.0173999965190887	-18922000\\
0.0179999992251396	-97694\\
0.0186000019311905	24682000\\
0.0191999971866608	47627000\\
0.0197999998927116	14204000\\
0.0204000025987625	939560\\
0.0209999978542328	-1970000\\
0.0216000005602837	-3773700\\
0.0222000032663345	-4634000\\
0.0227999985218048	-4512400\\
0.0234000012278557	-3454300\\
0.023999996483326	-2005200\\
0.0245999991893768	-277280\\
0.0252000018954277	1098400\\
0.025799997150898	1320700\\
0.0263999998569489	900770\\
0.0270000025629997	471610\\
0.02759999781847	97397\\
0.0282000005245209	-229790\\
0.0288000032305717	-249810\\
0.029399998486042	-166050\\
};
\addlegendentry{block 2}

\addplot [color=blue, line width=2.0pt]
  table[row sep=crcr]{%
0	0\\
0.000599995255470276	3235100\\
0.00120000541210175	6291100\\
0.00180000066757202	9584600\\
0.0023999959230423	736570\\
0.00300000607967377	-5511700\\
0.00360000133514404	-7389500\\
0.00419999659061432	-8028300\\
0.00480000674724579	-6666800\\
0.00540000200271606	-5112800\\
0.00599999725818634	-3827900\\
0.00660000741481781	419690\\
0.00720000267028809	1009700\\
0.00779999792575836	1851600\\
0.00839999318122864	2039300\\
0.00900000333786011	652080\\
0.00959999859333038	-310480\\
0.0101999938488007	-1789000\\
0.0108000040054321	-1528500\\
0.0113999992609024	-828540\\
0.0119999945163727	-2294900\\
0.0126000046730042	-6320600\\
0.0131999999284744	-14055000\\
0.0137999951839447	-18559000\\
0.0144000053405762	-12965000\\
0.0150000005960464	-2098800\\
0.0155999958515167	18550000\\
0.0162000060081482	45942000\\
0.0168000012636185	72164000\\
0.0173999965190887	-21819000\\
0.0180000066757202	-17515000\\
0.0186000019311905	-7931600\\
0.0191999971866608	2923700\\
0.0198000073432922	9314300\\
0.0204000025987625	6729100\\
0.0209999978542328	2696700\\
0.0215999931097031	714770\\
0.0222000032663345	-174950\\
0.0227999985218048	-612540\\
0.0233999937772751	-580480\\
0.0240000039339066	-387840\\
0.0245999991893768	-101010\\
0.0251999944448471	88411\\
0.0258000046014786	88051\\
0.0263999998569489	28410\\
0.0269999951124191	5708.59999999404\\
0.0276000052690506	3378.70000000298\\
0.0282000005245209	2963.70000000298\\
0.0287999957799911	1464.20000000298\\
0.0294000059366226	388.340000003576\\
};
\addlegendentry{block 5}

\addplot [color=red, line width=2.0pt]
  table[row sep=crcr]{%
0	0\\
0.000600000843405724	4199300\\
0.0011999998241663	7456600\\
0.00180000066757202	12582000\\
0.0023999996483326	13116000\\
0.00300000049173832	9373000\\
0.00359999947249889	1951900\\
0.00420000031590462	-5320000\\
0.00479999929666519	-12314000\\
0.00540000014007092	-15240000\\
0.00599999912083149	-12952000\\
0.00659999996423721	-6950300\\
0.00720000080764294	-1510300\\
0.00779999978840351	4990600\\
0.00840000063180923	8731800\\
0.00899999961256981	9231100\\
0.00960000045597553	8707200\\
0.0101999994367361	7219600\\
0.0108000002801418	6222000\\
0.0113999992609024	4861700\\
0.0120000001043081	3101900\\
0.0125999990850687	640200\\
0.0131999999284744	-1918800\\
0.0138000007718801	-3807100\\
0.0143999997526407	-5108600\\
0.0150000005960464	-6141100\\
0.015599999576807	-5963800\\
0.0162000004202127	-4937700\\
0.0167999994009733	-4317000\\
0.017400000244379	-5204400\\
0.0179999992251396	-6448300\\
0.0186000000685453	-6478300\\
0.0192000009119511	-4678700\\
0.0197999998927116	-2542000\\
0.0204000007361174	-3162600\\
0.0209999997168779	-2268300\\
0.0216000005602837	4355000\\
0.0221999995410442	3486000\\
0.02280000038445	3504300\\
0.0233999993652105	4756600\\
0.0240000002086163	113560\\
0.0245999991893768	-9428700\\
0.0252000000327826	-6454200\\
0.0258000008761883	1746400\\
0.0263999998569489	6374000\\
0.0270000007003546	834420\\
0.0275999996811152	-1502000\\
0.0282000005245209	-980070\\
0.0287999995052814	-248820\\
0.0294000003486872	-173610\\
};
\addlegendentry{block 4}

\addplot [color=mycolor2, line width=2.0pt]
  table[row sep=crcr]{%
0	0\\
0.000599998980760574	3686500\\
0.00120000168681145	6350100\\
0.00180000066757202	9239900\\
0.0023999996483326	3145200\\
0.00299999862909317	-3429400\\
0.00360000133514404	-5709300\\
0.00420000031590462	-6923000\\
0.00479999929666519	-6578000\\
0.00539999827742577	-4997100\\
0.00600000098347664	-3292600\\
0.00659999996423721	-713950\\
0.00719999894499779	617970\\
0.00780000165104866	1358900\\
0.00840000063180923	2409500\\
0.00959999859333038	355770\\
0.0102000012993813	-833630\\
0.0108000002801418	-171120\\
0.0113999992609024	1203700\\
0.011999998241663	3418400\\
0.0126000009477139	5184900\\
0.0131999999284744	6470500\\
0.013799998909235	6834600\\
0.0144000016152859	5639000\\
0.0150000005960464	1257600\\
0.0161999985575676	-14021000\\
0.0168000012636185	-17155000\\
0.017400000244379	-11720000\\
0.0179999992251396	2724800\\
0.0185999982059002	16800000\\
0.0192000009119511	23114000\\
0.0197999998927116	12826000\\
0.0203999988734722	4155200\\
0.0210000015795231	-557470\\
0.0216000005602837	-2251300\\
0.0221999995410442	-2919700\\
0.0227999985218048	-3025200\\
0.0234000012278557	-2263900\\
0.0240000002086163	-1559600\\
0.0245999991893768	-554740\\
0.0251999981701374	370990\\
0.0258000008761883	636500\\
0.0263999998569489	610060\\
0.0269999988377094	433240\\
0.0276000015437603	160740\\
0.0282000005245209	-37184\\
0.0287999995052814	-111620\\
0.029399998486042	-129390\\
};
\addlegendentry{block 3}

\end{axis}
\end{tikzpicture}%
		\caption{Pressure signal within the four connecting ``pillar''-like blocks 2, 3, 4, and 5 of the computational domain at time step 24400, showing the influence of different speeds of sound and damping parameters.\label{fig:Ex3_LinePlots}}
	\end{center}
\end{figure}

\paragraph{\textbf{Setup of the experiment}} Block 1 has a width and length of $0.025\,\textup{m}$ and a height of $0.01\,\textup{m}$; see \Cref{fig:Ex3_Domain}. On its bottom surface, the wave excitation takes place. We use an excitation signal in the form of a Dirichlet condition similarly to \Cref{subsec:Ex2}. The four walls of the block are equipped with homogeneous Neumann/symmetry boundary conditions, its top surface with a homogeneous Dirichlet condition, except at the interfaces to the blocks 2 to 5. Those are truncated four-sided pyramids which are aligned in a regular way between blocks 1 and 6, each with a height of $0.01\,\textup{m}$ and homogeneous Dirichlet conditions on its four walls. Block 6, with a height of $0.04\,\textup{m}$ and the remaining measures as for Block 1, covers the upper part of the geometry, again equipped with homogeneous Neumann/symmetry conditions on its four sides and homogeneous Dirichlet conditions at top and bottom, except for the interfaces with the blocks 2 to 5.\\
\indent The excitation signal is given by 
\begin{equation*}\label{eq:Excitation2}
g(t,x,y,z)=g^{(t)}(t)=\begin{cases}
\left(\frac{ft}{2}\right)^2\,A\,\sin(\omega t), & t<\frac{2}{f}\,\textup{s} \\
\hspace*{2.15cm}0, & t \geq \frac{2}{f}\,\textup{s};
\end{cases}
\end{equation*}
i.e., by a cut-off pulse-version of the continuous excitation signal used before. We use such a signal here in order to avoid interference within the block 1 originating from reflections off the walls between the interfaces to blocks 2 to 5. Amplitude and frequency are chosen as before. For the time discretization again the Generalized-$\alpha$-method is used with the final time $T=2.217\cdot 10^{-5}\,\textup{s}$, resolved with a step size of $dt = 10^{-9}\,\textup{s}$.\\
\indent \Cref{fig:Ex3_pics} shows snapshots of the solution computed with $p=2$ on 260730 elements. We observe iso-volumina of the highest acoustic pressure amplitudes at different time steps which show how the wave propagates through the four separate channels connecting the base block with the top block. Especially the deviations in the speed of sound are visible as the wave propagates much faster in block 4 (on the right) than, for example, in block 5 (on the left). This effect can also be seen in Figure \ref{fig:Ex3_LinePlots}, where the pressure signal is plotted along the central axes of the four ``pillars'' at a given time step.

While the signals corresponding to blocks 2 and 3 in Figure \ref{fig:Ex3_LinePlots} are traveling with the same speed (cf. Table \ref{tab:phys_param}), the signal in block 3 is much more damped compared to block 2, due to the damping parameter $b$ being much higher there. In contrast to that, the signal from block 5 is slower, while also higher in amplitude due to the changes in material properties and the signal in block 4 with the highest speed of sound has already passed through the ``pillar''-like structure and decayed in amplitude afterwards due to spreading into the empty space of block 6.

\begin{figure}[H]
	\begin{center}
		\includegraphics[trim=16cm 0mm 22cm 0mm, clip, scale=0.175]{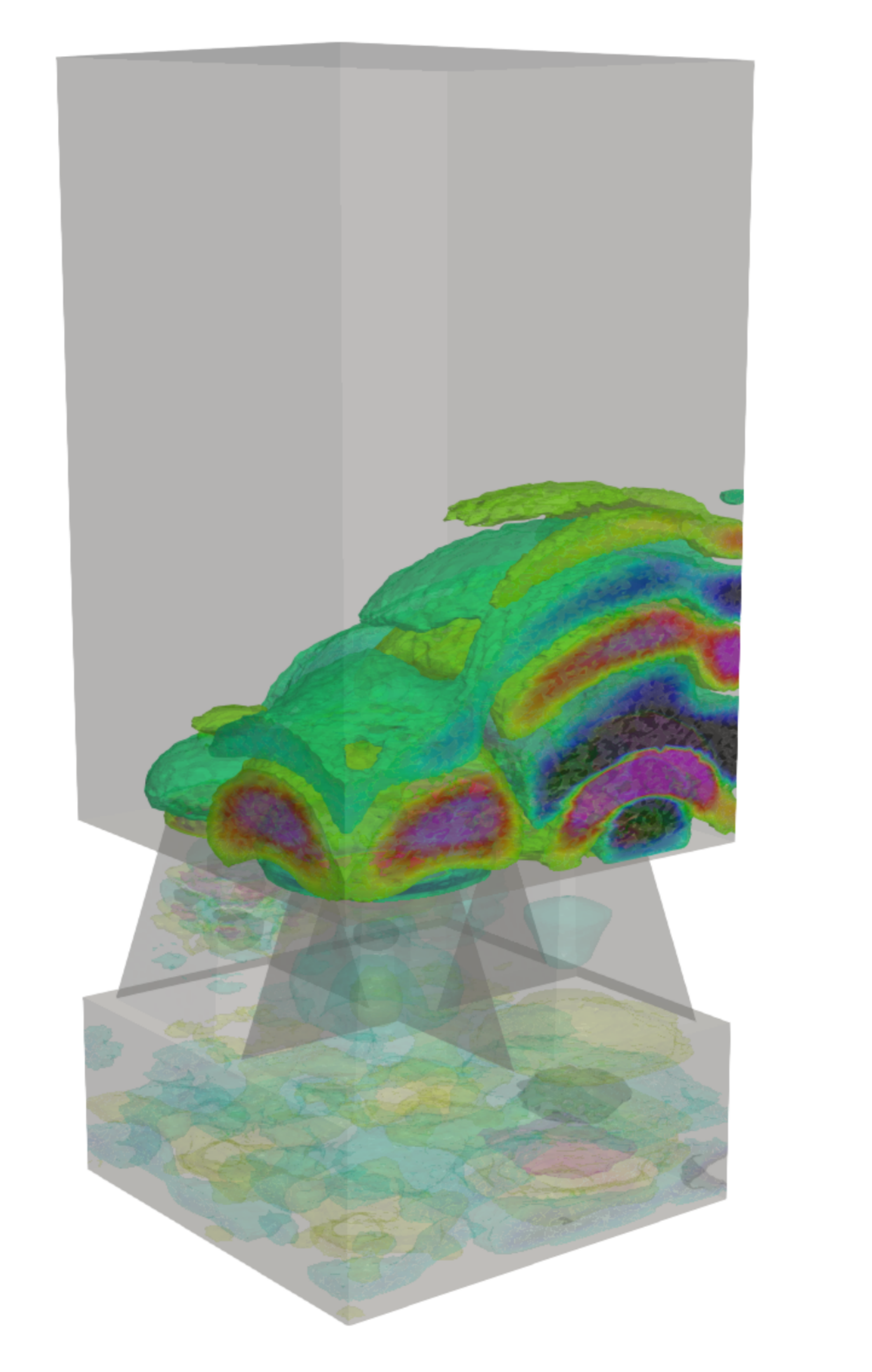}		\includegraphics[trim=16cm 0mm 22cm 0mm, clip, scale=0.175]{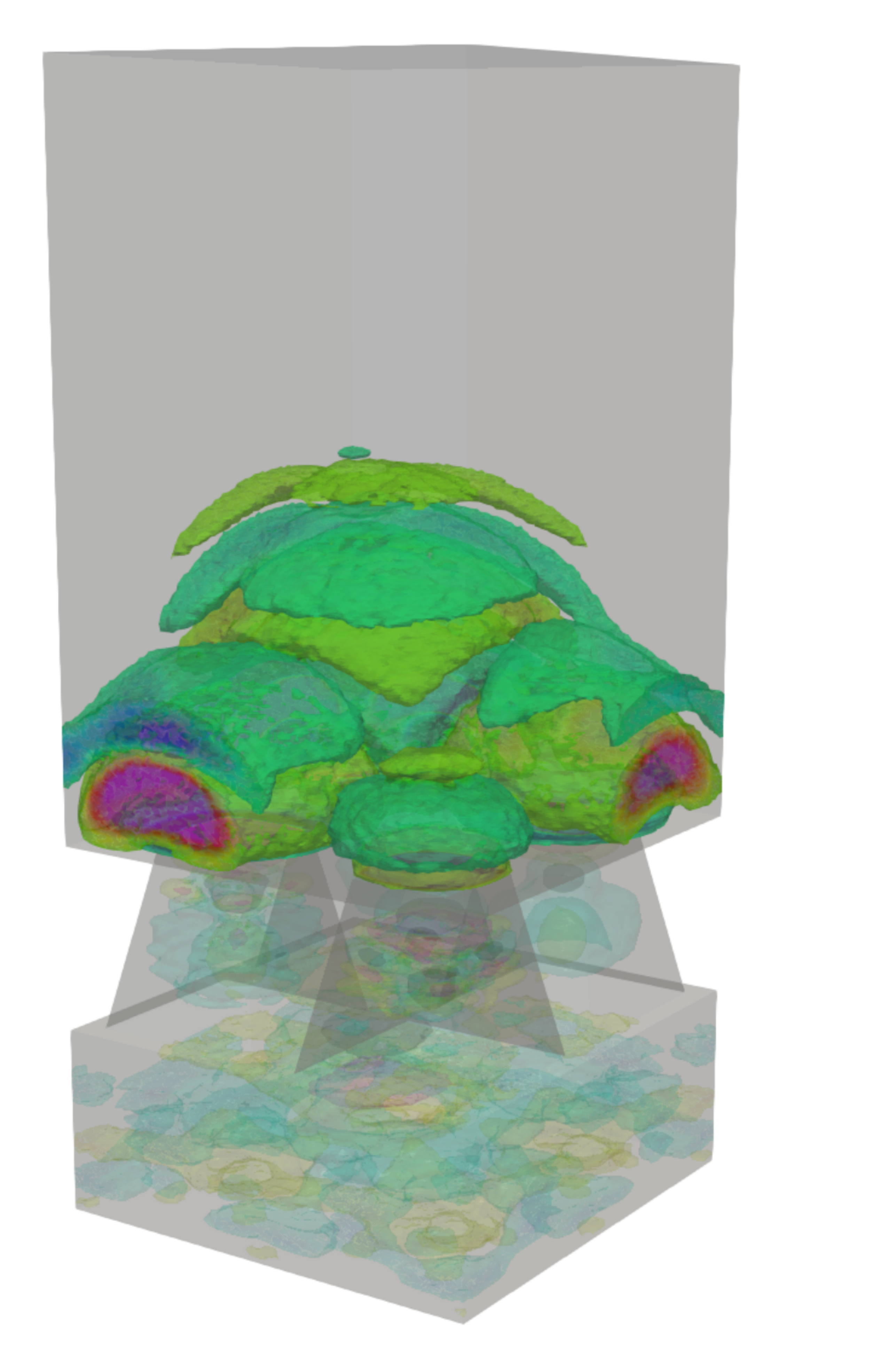}
		\raisebox{1.25cm}{\includegraphics[trim=0cm 0mm 0cm 0mm, clip, scale=0.3]{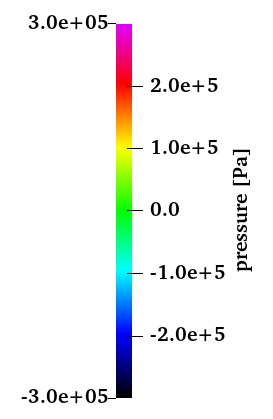}}
		\caption{Iso-volumina of the highest acoustic pressure amplitudes at time step 24600 within the upper part of the computational domain, where the individual waves coming from the connecting channels merge. \textbf{(left)} Same orientation of the domain as in \Cref{fig:Ex3_pics} \textbf{(right)} View rotated by $\approx \ang{90}$ at the same time step.\label{fig:Ex3_pics_farfield}}
	\end{center}
\end{figure}

\Cref{fig:Ex3_pics_farfield} shows iso-volumina of the acoustic pressure field as well, this time focusing on the wave propagation in the upper part of the computational domain, i.e., block 6 where the four individual waves coming from the four connecting channels again merge together into a single acoustic wave field. We observe decrease in the amplitude compared to \Cref{fig:Ex3_pics} due to the wave spreading.

\begin{figure}[H]
	\begin{center}
		\input{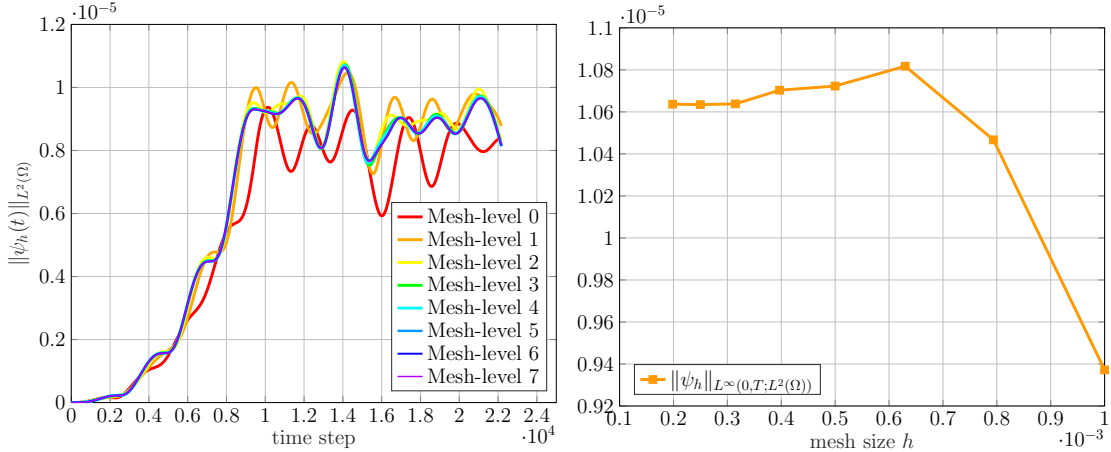}\hspace*{2mm}
%
%
\definecolor{mycolor1}{rgb}{1.00000,0.60000,0.00000}%
\begin{tikzpicture}[scale = 0.55, font=\Large]

\begin{axis}[%
width=4.568in,
height=3.589in,
at={(0.766in,0.5in)},
scale only axis,
xmin=0.0001,
xmax=0.001,
xlabel style={font=\Large\color{white!15!black}},
xlabel={mesh size $h$},
ymin=9.2e-06,
ymax=1.1e-05,
ylabel style={font=\Large\color{white!15!black}},
ylabel={},
axis background/.style={fill=white},
title style={font=\bfseries},
title={},
xmajorgrids,
ymajorgrids,
legend style={at={(0.03,0.03)}, anchor=south west, legend cell align=left, align=left, draw=white!15!black}
]
\addplot [color=mycolor1, line width=2.0pt, mark=square*, mark options={solid, fill=mycolor1, mycolor1}]
  table[row sep=crcr]{%
0.001	9.37207257967481e-06\\
0.0007937005259841	1.04669938542608e-05\\
0.000629960524947437	1.08175052077107e-05\\
0.0005	1.07230230131363e-05\\
0.00039685026299205	1.07032777211628e-05\\
0.000314980262473718	1.06382186876896e-05\\
0.00025	1.06348204587825e-05\\
0.000198425131496025	1.06365745569789e-05\\
};
\addlegendentry{$\|\psi_h\|_{L^{\infty}(0,T;L^2(\Omega))}$}

\end{axis}
\end{tikzpicture}%
		\caption{Quantity of interest over an $h$-refinement with quadratic shape functions. \textbf{(left)} $L^2(\Omega)$ norm of the numerical solution over the course of time for all tested levels of refinement \textbf{(right)} Respective $L^{\infty}(0,T;L^2(\Omega))$-norm results on all 8 mesh levels in order to observe convergence. \label{fig:Ex3_href}}
	\end{center}
\end{figure}

As before, exact solution is not available. Therefore, we again track the quantity of interest given by \eqref{QoI}.  We evaluate it over an $h$-refinement with quadratic shape functions on meshes with mesh-sizes of $h_j = 0.001\sqrt[3]{2^{-j}},~j=0,...,7$. The results are depicted in \Cref{fig:Ex3_href}.

We observe a convergence of the quantity of interest towards a value of around $1.064\cdot 10^{-5}$ as $h$ approaches zero; see~\Cref{fig:Ex3_href} on the right. We note, however, that allowing for the jumping material coefficients lies beyond the theory presented in this work. Therefore, we can conclude that the application of the spectral discontinuous Galerkin method on a problem with varying coefficients is feasible, while a rigorous convergence analysis is left for future work.

%
\section*{Acknowledgements}
P. F. Antonietti has been partially funded by PRIN grant n. 201744KLJL funded by the  Ministry of Education, Universities and Research (MIUR). Paola F. Antonietti and I. Mazzieri  also acknowledge the financial support given by INdAm-GNCS. \\ \indent M. Muhr, V. Nikoli\'c, and B. Wohlmuth acknowledge the financial support  provided  by the Deutsche Forschungsgemeinschaft under the grant number WO 671/11-1.
\bibliographystyle{abbrv} 
\bibliography{references}
\end{document}